\documentclass[a4paper,11pt,english,twoside, fleqno]{amsart}

		\usepackage{amsfonts, amsmath, amssymb}
	\usepackage[fixamsmath]{mathtools} 
	\usepackage{stmaryrd} 
	\usepackage[all]{xy}
	\usepackage{graphicx}

	\usepackage[utf8]{inputenc}
	\usepackage[T1]{fontenc}
	\usepackage{lmodern}
	\usepackage{cmap} 

	\usepackage[textwidth=15cm,heightrounded,hcentering]{geometry}
	\usepackage{wrapfig}
	\usepackage{calc}
	\usepackage[usenames,dvipsnames,svgnames,table]{xcolor}
	\usepackage{color}
	\definecolor{e-mail}{rgb}{0,.40,.80}
	\definecolor{reference}{rgb}{.20,.60,.22}
	\definecolor{mrnumber}{rgb}{.80,.40,0}
	\definecolor{citation}{rgb}{0,.40,.80}
	\definecolor{gris25}{gray}{0.60}

\usepackage{xspace}
\usepackage[protrusion=true,expansion=true,kerning, spacing]{microtype}
\microtypecontext{spacing=nonfrench}
\usepackage[english]{babel}

	\theoremstyle{plain}
	\newtheorem{lemma}{Lemma}[subsection]
	\newtheorem{theorem}[lemma]{Theorem}
	\newtheorem{theoNoSub}{Theorem}[section] 
	\newtheorem{proposition}[lemma]{Proposition}
	\newtheorem{remark}[lemma]{Remark}
	\newtheorem{definition}[lemma]{Definition}
	\newtheorem{defiNoSub}[theoNoSub]{Definition} 
	\newtheorem{example}[lemma]{Example}

	\newcommand{\Mg}{\mathcal{M}_g}
	\usepackage{xparse}
	\DeclareDocumentCommand{\Hurw}{O{G} O{g}} {\mathcal{M}_{#2}[#1]}
	\DeclareDocumentCommand{\Mgm}{ O{g} O{m}}{\mathcal{M}_{#1,#2}}
	\DeclareDocumentCommand{\Mzero}{ O{m}}{\M[0][#1]}
	\DeclareDocumentCommand{\Mdmgm}{ O{g} O{m}}{\overline{\mathcal{M}}_{#1,#2}}
	\DeclareDocumentCommand{\Mdmzero}{ O{m}}{\Mdm[0][#1]}
	
	\DeclareDocumentCommand{\res}{O{G} O{H}} {\mathrm{res}^{#1}_{#2}}
	\DeclareDocumentCommand{\cores}{O{G} O{P}} {\mathrm{cores}^{#1}_{#2}}
	
	\newcommand{\QQ}{\mathbb Q}
	\newcommand{\ZZ}{\mathbb Z}
	\newcommand{\Gq}{\textrm{Gal}(\bar{\QQ}/\QQ)}
	\newcommand{\M}{\mathcal M}
	\newcommand{\spec}{{\mathrm{Spec}\,}}
	\newcommand{\aut}{{\mathrm{Aut}}}
	
	\usepackage{leftidx}
	\newcommand{\Diam}{\lozenge}
	
	\usepackage{tikz-cd}

	\usepackage{stmaryrd}
	\usepackage[sans]{dsfont} 
	\usepackage{bm} 
	\usepackage{booktabs} 
	\usepackage{enumitem}
	\setenumerate[0]{label=\roman*)}

	\definecolor{gris25}{gray}{0.60}\definecolor{Gray85}{gray}{0.85}\definecolor{Gray65}{gray}{0.65}
	
	\usepackage{lastpage}
	\usepackage{fancyhdr}
	\newcommand{\et}{{\textrm{\scriptsize\'et}}}
	\renewcommand{\H}{{\mathrm{H}}}
	\newcommand{\R}{{\mathrm{R}}}
	\newcommand{\dquot}{/\hskip-0.2em/}

	\usepackage{titletoc}
	\usepackage[explicit]{titlesec}
	\newcommand\dotafter[1]{\if\relax\detokenize{#1}\relax{}\else {#1.}\fi}
	
	\titleformat{name=\section}[block]{\centering\large\scshape}{}{0em}{\vspace*{0em}}[\thesection.~#1]
	\titleformat{name=\section,numberless}[block]{\centering\large\scshape}{}{0em}{\vspace*{.2em}#1}[] 
	\titleformat{\subsection}[runin]{\bfseries}{\thesubsection.~}{0em}{}[\dotafter{#1}]
	\titlespacing{\subsection}{0pt}{*2}{1ex}
	\titleformat{\subsubsection}[runin]{\itshape}{\thesubsubsection.~}{0em}{}[\dotafter{#1}] 

	\titlecontents{section}[.5em]{\vspace{.3em}\bfseries\fontfamily{phv}\selectfont}{\thecontentslabel.~}{}{\titlerule*[1pc]{}\contentspage}
	\titlecontents{subsection}[3em]{}{\contentslabel{.7cm}}{\hspace*{2.3em}}{\titlerule*[1pc]{.}\contentspage}
	
	\usepackage{standalone}
	
	\usepackage[pdftex,%
	backref, 
	hyperindex=true, bookmarks=true, bookmarksopen=true, bookmarksnumbered=true,
	pdfauthor={B. Collas, S. Maugeais},
	pdftitle={Hurwitz Stacks and Irreducibility},
	pdfsubject={},
	pdfkeywords={},
	pdfstartview=FitH,
	colorlinks=true,
	linkcolor=reference, citecolor=citation, urlcolor=e-mail
	]{hyperref} 

	\usepackage{caption}
\usetikzlibrary{fit,backgrounds} 

	\usepackage[alphabetic,initials,lite]{amsrefs}
	
	\def\MR#1{ \href{http://www.ams.org/mathscinet-getitem?mr=#1}{MR#1}}

	\numberwithin{equation}{subsection} 

	\title[Irreducibility of Hurwitz Stacks]{Hurwitz Stacks of Groups Extensions and Irreducibility}
	\author[B.~Collas, S.~Maugeais]{Benjamin Collas, Sylvain Maugeais}
	\keywords{Hurwitz stack, moduli stack of curves, special loci, mixed cohomology, deformation of curves}

	\subjclass[2010]{ %
		14H10, 
		14H30 (primary), and 
		14H45, 
		32G15 
		(secondary).} 
	
	\thanks{This research was partially supported by DFG programme DE 1442/5-1, DFG programme SPP 1786 and FR CNRS 2962, and by the Research  Institute  for  Mathematical  Sciences of Kyoto University . The first author also gratefully acknowledges an invitation to the programme \emph{Algebro-Geometric and Homotopical Methods} at the Mittag-Leffler Institut, as well as M.~Fried for some precious comments on a first version of this work.
}

\begin{document}
	\addtocontents{toc}{\protect\enlargethispage{2cm}}

	\begin{abstract}
	We study the irreducible components of special loci of curves whose group of symmetries is given as certain group extension. We introduce some relative Hurwitz data, which we show by using mixed étale cohomology theory, identifies some irreducible components for rational and normal non-abelian special loci and Hurwitz spaces. A heuristic, that is supported by three classes of examples, provides an additional context for building further irreducible loci.
	\end{abstract}

	\vspace*{-2em}
	\maketitle
	
	\vspace*{-2em}

	\tableofcontents
	\thispagestyle{empty}
	\newpage

\section{Introduction}
	 This paper is motivated by the study of the \emph{stack arithmetic of the moduli spaces of curves} of genus $g$ with $m$-marked points $\M_{g,[m]}$, i.e the study of the $\Gq$-action on its stack inertia groups $I_{\M,x}(\bar \QQ)$, and more precisely by its geometric formulation in terms of \emph{irreducible components of special loci $\M_{g,[m]}(G)$} of curves with automorphism group $G\simeq I_{\M,x}(\bar \QQ)$. As detailed in \cite{ColMau2} \S 1.1, this approach relies on determining \emph{algebraic invariants in family of the irreducible components of $\M_{g,[m]}(G)$}. Via the identification of the normalisation $\tilde{\M}_g(G)\simeq\Hurw/\aut(G)$, this question can be reformulated at the level of Hurwitz spaces $\Hurw$ of curves endowed with a $G$-action, and thus shares some similarity, in motivations and techniques, with the role of $G$-covers in the \emph{Regular Inverse Galois Problem}.
	
	\medskip
	
	After a brief reminder on the arithmetic of Hurwitz spaces and special loci, we first discuss the definition of our relative invariants $\Diam$ of components of G-covers with rational locus  $\M_g[G]^{rat}$, then present our irreducibility result for \emph{rational special loci} $\M_g(G)^{rat}_\Diam$, and finally introduce certain Heuristic that is supported by $3$ types of loci and produces multiple irreducible special loci $\M_g(G)^{rat}_\Diam$ and $\M_g(G)_{\bar \Diam}$.

	\subsection{Hurwitz Spaces and Special Loci in Arithmetic Geometry} \label{sub:IntroLoci}
		Let $G$ be a finite group, and denote by $\Hurw$ and $\Mg(G)$ respectively the Hurwitz stack of $G$-covers and the special loci of $G$ in $\Mg$ -- i.e. the moduli stack of curves admitting a $G$-action. Arithmetic motivations for the study of the irreducible components of those stacks comes from two similar questions in the study of the absolute Galois group of rational numbers $\Gq$.
		
		\medskip

		The \emph{Regular Inverse Galois Problem} (RIGT) -- that is to realize $G$ as a quotient of the absolute Galois group of rationals $\Gq$ -- turns into the geometric questions of realizing $G$ as a regular Galois cover of $\mathbb{P}^1\setminus\{0,1,\infty\}$ and to the existence of a $\QQ$-rational point in $\Hurw$. Weakening the latter Diophantine condition to the question of the existence of \emph{geometrically irreducible components in $\Hurw$} leads to the construction of Harbater-Mumford irreducible components that are defined over $\QQ$ -- see \cite{FRI95} --, then the realization of every center-free group $G$ over $\mathbb F_q$ -- see \cite{WEW98}, and also \cite{DD06} for the realization of a projective system of finite groups in terms of Hurwitz towers.

		\medskip

		The \emph{Geometric Galois Actions} (GGA) study the actions of $\Gq$ on the étale fundamental group of moduli stacks of curves $\Gq\to Aut[\pi_1^{et}(\M_{g,[m]}\otimes \bar{\QQ}),\bar x]$, where $\bar x$ denotes a geometric point of $ \M_{g,[m]}$. The question of the $G$-arithmetic appears through the stack inertia $I_{\M,\bar x}(\bar \QQ)\simeq G\hookrightarrow \pi_1^{et}(\M_{g,[m]}\otimes \bar{\QQ},\bar x)$ or equivalently the irreducible components of special loci $\M_{g,[m]}(G)$ -- see \cite{ColMau2} \S 1 for details. This geometric approach leads to the description of the $\Gq$-action for respectively $p$-groups and cyclic groups \cite{ColMau1} Theorem 5.9, \cite{ColMau2} Theorem 4.8, where it is proven to be given by $\chi$-conjugation.
		
		The general context being given by the inertia stratification in local gerbes bounded by the automorphism group of objects, this raises the question of \emph{describing this $\Gq$-action on higher inertia groups}. A first essential step in this direction is given by the \emph{definition of rational invariant of irreducible components of the special loci in family} -- see \cite{ColMau1} Theorem 4.3 for the cyclic case and for quotient curves of any genus.
		
		\medskip
		
		In this paper, we exploit the arithmetic of $G$-covers through the rationality of the ramification locus:  We deal with the stack $\M_{g}(G)^{rat}_{\Diam}$ of \emph{rational special loci} with rational ramification locus and with given rational Hurwitz datas $\Diam$ that extends the invariants $\bf k$ built in \cite{ColMau1} \S 3.1 for the need of the cyclic study  -- see Remark~\ref{rem:HurwCaseCyc}.

	\subsection{Invariants of Components of Special Loci}\label{sub:IntroArith}
		The construction of discrete invariants of irreducible components of $\M_{g}(G)$ relies on a long tradition of \emph{geometric invariants of $G$-covers} via group theoretic methods in terms either of equivalence classes of \emph{generating vectors} -- see \cite{BRO90a} for $G$-curves of genus $2$ and $3$, and \cite{BROWO07} for abelian groups --, or of equivalence classes of monodromy representations of the $G$-quotient curve $D$ as for RIGT.

		Such numerical invariants include: (1) the genus of $D$, (2) the Nielsen invariants counting the number of local monodromy belonging to a given set of conjugacy classes, and more recently (3) a \emph{global second homology invariant} \cite{CLP15} which refines the previous ones in the case of $G$-action with étale factorization and $G=D_{2n}$ -- see ibid.

		\medskip
		
		Our approach is to define some similar \emph{relative} invariants $\bf k$ and $\Diam$, i.e. for a family $C/S$ of $G$-curves in $\M_{g}[G]$. For cyclic group, the construction of $\bf k$ is achieved in terms of étale cohomology, see \cite{ColMau1} \S 3, with similar irreduciblity results for $\M_g(\gamma)_{\bf k}$ -- see Theorem 4.3 and Proposition 3.12 ibid. Their generalisation $\Diam$ and $\bar\Diam$ in \S \ref{subsub:RatHurwSt} involves some additional characters data $\{\chi_i\}_I$ and supports some \emph{arithmetic properties} -- see \S \ref{subsub:GArtih}.
	
		The irreduciblity property is then dealt with by the use of mixed cohomology $H_P^{\bullet}(X_{et},H)$ of \cite{ClassesChern}[Chap. V] in the case of $(H;G)$-covers, $G$ being an extension of $P$ by $H$: this determines both the ramification and the global parts of the covers within an algebraic deformation -- see \S \ref{sub:DefLoc} and \S \ref{sub:DefGlob} for the construction of this deformation. Under some additional assumptions, this in turns provides some irreduciblity results for $\M_g(H\triangleleft G)^{rat}_\Diam$ and $\M_g(H \triangleleft G)_{\bar \Diam}$

	\subsection{Characterizing Irreducible Components of Special Loci}\label{sub:IntroIrrSpecLoci}
		Under a certain $\aut_\Diam$-liftability condition of $G$ with respect to $\Diam$ and $H\triangleleft G$, the main result of this paper is to spread the irreducibilty from $\M_{g'}(G/H)^{rat}_{\textrm{cores}(\Diam)}$ to $\M_g(G)^{rat}_{\Diam}=\mathcal M_g[G]^{rat}_{\Diam}/\aut_{\Diam}(G)$ -- see Theorem~\ref{th:IrredGeneral} and above for the precise definitions:
 
		\newtheorem*{theorem*}{Theorem}
		\begin{theorem*}[A]
			Let $G$ be a finite group and $\Diam$ a $G$-Hurwitz data so that $G$ is $\aut_{\Diam}$-liftable with respect to $H\simeq\ZZ/p\ZZ\triangleleft G$, that $\Diam$ is $H$-irreducible, and that $\Diam$ contains $H<G$ as isotropy group.  Then $\M_g(G)^{rat}_{\Diam}$ is irreducible.
		\end{theorem*}
	
		This result relies essentially on the irreducibility of the relative quotient morphism $\Phi_\Diam^{rat}\colon \mathcal M_g[G]^{rat}_{\Diam}/\aut_{G/H, \Diam}(G;H) \to \mathcal M_{g'}[G/H]^{rat}_{\cores(\Diam)}$, see Theorem \ref{th:RelativeIrr}. Under the existence of a $H$-ramification point on top of a $P$-ramification one, we proceed by building a certain deformation via the compactification of $G$-equivariant torsors of \cite{MaugeaisBordeaux}, see Proposition \ref{prop:DefBrLoci}, and the identification of $(H;G)$-covers to a certain stack $[R^1f_*H]^{G/H}$ of $H$-torsors, see \S \ref{subsub:CompGTors}. Notice that the deformation of the global part relies on an analytic result that may be unavoidable since it is already a key ingredient in proving the irreducibility of $\Mg$ in \cite{DM69} and of the special loci $\M_{g,[m]}[\ZZ/n\ZZ]_{\lozenge}$ in \cite{COR87,CAT12} then \cite{ColMau1}.
	
		\medskip
		
		We furthermore introduce \emph{an $\aut_{\Diam}$-solvable Heuristic} -- that is the existence of a series $G = G_0 \triangleright G_1 \triangleright \cdots \triangleright G_m = \{1\}$ that provides a sequence of data $\{(\lozenge_k,H_k)\}_{1\leqslant k \leqslant m}$ whose associated intermediate quotients satisfy the hypothesis of Theorem (A), see \S \ref{subsub:AutLozHeur} and Fig.~\ref{fig:ItIrr} -- which provides further irreducible rational special loci, see loc. cit. and Remark~\ref{rem:MoreIrrLoci}. We support this Heuristic by introducing $3$ types of special loci: (1) a \emph{cyclic type} $\bf C$, (2) an \emph{Elementary Abelian type $\bf{AE}$}, and (3) a \emph{dihedral type} $\bf D^{rot}$ -- see \S \ref{subsub:CorAbDi}. Our approach partially covers the dihedral cases of \cite{CLP15} -- see Remark~\ref{rem:D2nIrr}, and is complementary to group cohomology methods (e.g. Schur multiplier and thus Fried's universal Frattini covers) -- see Remark~\ref{rem:GrpCohomIsNoGood}.
		
		Moreover, we discuss the question of the descent from $\M_g(G)_\Diam^{rat}$ to $\M_g(G)_\Diam$ that tends to show that the arithmetic invariants $\Diam$ provide either a rational irreduciblity for $\M_g(G)_\Diam^{rat}$, or descend to $\M_g(G)_{\bar \Diam}$ -- see \S \ref{subsub:AutLozHeur} and the $D_{2p}$-case in \S \ref{sub:TypeDrot}.
	
		\bigskip

		The plan of the article goes as follow. The definition of rational Hurwitz data and rational Hurwitz loci, as well as moduli properties by (co)restriction and under $\aut$-quotient are presented in \S 2; the development of mixed cohomology and their relation with $(H;G)$-covers and torsors occupies \S 3; the relative irreducible result for $\Phi_\Diam^{rat}$ is established in \S 4 following the deformation of $(H;G)$-covers with $H$-ramification -- with algebraic and analytical results; \S 5 deals with the irreducibility of the rational special loci $\M_g(G)^{rat}_\Diam$, also discusses the $\aut_{\Diam}$-liftability Heuristic, and illustrates the non-triviality of the $\aut$-liftability property; \S 6 presents the various types that support the $\aut_{\Diam}$-liftability Heuristic. An Appendix corrects a mistake of taking a stack quotient in \cite{MaugeaisBordeaux}.

		\bigskip
		
		We expect this paper, in combination with a thoughtful study of explicit groups and Hurwitz data and with the use of inertial limit Galois actions of \cite{ColMau2}, to lead to further descriptions of the GGA on \emph{higher stack inertias} $I_{\M,x}(\bar \QQ)$. 
	
\section{Hurwitz Characters and Normality}
	Let $G$ be a finite group and $g \ge 0$ an integer. We consider the stack $\Hurw$ of smooth proper curves of genus $g$ with a given faithful $G$-action, which is a Deligne-Mumford stack over $\ZZ[1/|G|]$. {In what follows, the scheme $S$ belongs to the category of $\ZZ [1/|G|]$-schemes or $\QQ$-schemes to avoid any wild ramification difficulties.}

	\subsection{Rational Hurwitz Invariants} 
		Let $C/S$ be a smooth proper curve endowed with a faithful $G$-action $\iota\colon G\hookrightarrow \aut_S(C)$, i.e. faithful in every fibre. After some reminders on the reduced ramification locus $Ram_S(C,\iota)$ of $(C/S,\iota)$, we define a notion of \emph{rational Hurwitz data} over $S$ which are shown to reflect rational and arithmetic properties of $(C/S,\iota)$ and is a variant of the geometric ones introduced in \cite{FRIED77} and \cite{SER08}.

		\subsubsection{} \label{subsub:ReducedRam}
			Let	$\mu_{|G|}$ denote the group scheme of $|G|$-roots of unity, and consider a geometric fibre $C_{\bar s}$ of $C/S$ at $\bar s=\spec (k)$. For every ramified point $x\in C_{\bar s}$, denote by $G_x<G$ its isotropy group. This is a cyclic group since the action is tame, and it defines a \emph{Hurwitz character} via the cotangent representation:
			 \[
			  \chi_x\colon G_x\to \mu_{|G|}(k)
			 \]  
			 that is locally given by  $\chi_x(h)=h(\varpi)/\varpi\mod \varpi^2$, where $\varpi$ is a uniformising parameter of $C$ at $\bar s$. This representation is primitive, i.e. of order $|G_x|$. 
			
			\bigskip
			
			In the relative case, denoting $C^{H}$ the scheme of fixed points, one defines a ramification divisor
			\[
			Ram(C,\iota)=\sum\nolimits_{H\neq 1}\phi(|H|) C^{H} \text{ where } \phi \text{ denotes the Euler function}
			\]
			which is a \emph{relative Cartier divisor} over $S$ -- we refer to \cite{BERO07} \S 3 for the original approach and details.
			
			Suppose that the support of the ramification divisor $Ram(C,\iota)$ is given by disjoint sections $(e_i)_{i \in I}$. By the same arguments as above, one attaches to each $e_i$ a cyclic stabilizer group $G_i$ and a \emph{primitive character}
				\begin{equation}\label{eq:Character}
					\chi_i \colon G_i \to  \mu_{|G|},
				\end{equation} 
			which are locally constants over $S$ -- see especially Lemme 3.1.3 Ibid. If $S$ is the spectrum of an algebraically closed field, these characters are induced by the action of $G_i$ on the tangent space of $C$ at $C(e_i)$.

			\bigskip

			Let us introduce the stack $\Mgm[g][I][G]^{rat}$ of $G$-covers with rational ramification locus.

			\begin{definition}
				Let $G$ be a finite group and $I$ a finite $G$-set of cardinal $m$. We denote by $\Mgm[g][I][G]^{rat}$ the stack classifying equivariant $G$-curves $(C/S,\iota)$ together with disjoint sections $\{e_i\colon S \to C\}_{i=1,\dots,m}$.
			\end{definition}
			
			For such a $G$-curve, one has 
			\begin{equation}\label{eq:Ram}
				Ram(C,\iota) = \bigoplus_{i=1}^m \phi(|G_i|) e_i,
			\end{equation}
			where $G_i$ is the stabiliser of $e_i$ in $G$.

			\medskip
			
			Denoting $|I|=m$, the stack $\Mgm[g][I][G]^{rat}$ is naturally a closed substack of the algebraic stack of $m$-pointed $G$-curves $\Mgm[g][m][G]$. Moreover, as the ramification divisor of every $G$-curve $(C/S,\iota)$ splits up to an étale base change, the morphism given by forgetting the ramification divisor
			\[
			\Mgm[g][I][G]^{rat} \to \Mg[G]
			\]
			is étale.
			
			The equation \eqref{eq:Ram} tells us in particular that the genus $g'$ of the quotient curve $C/G$ is fixed, as it is given by the Hurwitz formula:
			\[
			2g-2 = (2g'-2)|G| + \sum_{i \in I} (|G_i| - 1)
			\]

			In a similar way to the moduli stack of curves, one proves that the stack $\Mgm[g][I][G]^{rat}$ is {a smooth Deligne-Mumford algebraic stack as GIT quotient} -- see \cite{BERO07} \S 6.3. We recall that $\Hurw$ (resp.  $\Mgm[g][I][G]^{rat}$) is not a substacks of $\Mg$ (resp. of $\Mgm[g][m]$), and that forgetting the $G$-action corresponds to the normalization morphism $\Hurw/\aut(G)\to \Mg(G)$ (resp. idem for $m$-pointed curves) -- see \cite{Rom09} \S 3.4 and \cite{ColMau1} \S 2.1 for details.

			\medskip
			
			The notion of \emph{Hurwitz characters} encodes some geometric and Galois properties that we formalize in the next section. 

		\subsubsection{} \label{subsub:RatHurwSt}
			Let $(C/S,\iota)\in \Mgm[g][I][G]$ be an equivariant curve, with isotropy groups $\{G_i\}_I$, to which can be associated some primitive characters $\{\chi_i\}_I$ in a functorial way as in Eq.~\eqref{eq:Character} above. {We recall that for a $\gamma \in G$, the $G$-action on $C$ induces conjugacy morphisms $\psi_{\gamma,i}\colon G_i\to G_{\gamma.i}$ on the isotropy groups which on the characters translates as follow:
				\begin{equation}\label{eq:HurwCharRel}
					G_{\gamma.i}=\gamma^{-1}G_i\gamma\text{ and }\chi_{\gamma.i}=\chi_i\circ \psi_{\gamma,i}^{-1}\text{ where } \gamma\in G.
				\end{equation}
			}
			
			The following is a variant of \cite{SER08, FRIED77} which originally defines the Hurwitz data on geometric fibre and as conjugacy classes of Hurwitz datas.

			\begin{definition}
				Let $G$ be a finite group and $I$ a finite $G$-set of cardinal $m$. A \emph{Hurwitz data} $\Diam$ associated to $G$ and $I$ is a $m$-tuple $((G_1,\chi_1),...,(G_{m},\chi_{m}))$ of cyclic subgroups $G_i<G$ and primitive $G_i$-characters $\chi_i\in \mathrm{Isom}(G_i,\mu_{|G|})$ such that equations \eqref{eq:HurwCharRel} are satisfied. 
			\end{definition}
			
			We denote by $\Diam{_G^I}$ the set of $\Diam{_G^I}$-Hurwitz data associated to $G$ and $I$. For $\Diam\in\Diam_G^I$, we write $stab(\Diam)$ (resp. $char(\Diam)$) for the list of $G$-isotropy groups (resp. $G$-primitive character), {and we call $I$ the set of indices of $\Diam$.}

			The Hurwitz data attached to a $G$-cover of $\Mgm[g][I][G]^{rat}$ encodes all its classical elementary invariants (eg. the ramification indices and the genus of the quotient curve $C/G$). For a cyclic group $G$, it also coincides with the $\gamma$-type introduced in \cite{ColMau2} \S 2.2 whose relation to the classical branching data $\{k_j\}_I$ of a geometric cover is given by Lemma 2.4 Ibid -- see also Remark~\ref{rem:HurwCaseCyc}.

			\medskip
			
			Let $\Diam$ be a $\Diam{_G^I}$-Hurwitz data, we denote by $\Mgm[g][I][G]_{\Diam}^{rat}$ the substack of $\Mgm[g][I][G]^{rat}$ whose sections are $G$-curves with Hurwitz datas $\Diam$. 

			\begin{proposition}
				The stack $\Mgm[g][I][G]_{\Diam}^{rat}$ is open and union of connected components in $\Mgm[g][I][G]^{rat}$.
			\end{proposition}
			
			This immediately follows first from the local constance of the Hurwitz data by \cite{BERO07} Lemme 3.1.3, then from the stability under specialisation.
			
			\medskip
			
			{Recall that the set $\Diam{_G^I}$ is naturally endowed with a $G$-action through the action of $G$ on $I$ induced by the conjugacy action of Eq. \eqref{eq:HurwCharRel}.}

			\begin{definition} 
				A Hurwitz data $\Diam = ((G_1,\chi_1),...,(G_{m},\chi_{m}))$ is said to be \emph{normal} if for all $i \in I$ and $\gamma \in G$ we have $G_{\gamma.i} = G_i$ and $\chi_{\gamma.i} = \chi_i$.
			\end{definition}
			
			The choice of a representative in each class in $\Diam{_G^I}/G$ is then sufficient to recover the whole Hurwitz data \emph{up to an order on $I$ only}. For a normal Hurwitz data $\Diam$ and for $\gamma \in G$, the difference between $\Diam$ and $\gamma.\Diam$ is only the permutation of $I$ given by $\gamma$.
			
			\begin{proposition}\label{Prop:normAb}
				Let $\Diam = ((G_1, \chi_1),\dots,(G_m,\chi_m))$ be a $\Diam{_G^I}$-Hurwitz data. If $G_i < Z(G)$ for all $i\in I$, then $\Diam$ is normal. In particular, every Hurwitz data of an abelian group is normal.
			\end{proposition}
			
			\begin{proof}
				It follows from Eq.\eqref{eq:HurwCharRel} that a Hurwitz data $\Diam$ is normal if and only if for all $i\in I$ one has:
				\[
				\chi_i=\chi_{\gamma.i},\, \forall \gamma \in G\ \text{ i.e. }\
				\begin{cases}
				G_i\unlhd G\\
				\chi_{i}(h)=\chi_i(\gamma h \gamma^{-1}),\forall h\in G_i,\, \forall \gamma \in G,
				\end{cases}			 
				\]
				which by the injectivity of the characters, is equivalent to $G_i< Z(G)$ for all $i\in I$.				
			\end{proof}
			
			\medskip

			In the case of curves in $\M_g[G]$ with non-rational ramification locus, an \emph{unordered Hurwitz data} is given as a class $\bar\Diam \in \Diam{_G^I}/G$, and we define $\M_{g}[G]_{\bar\Diam}$ as the locus of $G$-curves in $\M_{g}[G]$ that are étale locally isomorphic to a curve inside $\Mgm[g][I][G]_{\Diam}^{rat}$, for a certain $\Diam \in \bar \Diam$. The stack $\M_g[G]_{\bar\Diam}$ is naturally a Deligne-Mumford algebraic substack of the curves with $m$-marked points $\M_{g,[m]}$. As the ramification divisor splits \'etale locally in a $G$-equivariant way, we have an \'etale surjective morphism 
			\[
			\coprod_{\Diam \in \bar \Diam} \Mgm[g][I][G]_{\Diam}^{rat} \to \M_g[G]_{\bar\Diam}.
			\]

			\medskip
		
			We have the following proposition.

			\begin{proposition}\label{prop:OrdUnordNorm}
				Let $G$ be a finite groupe, $I$ a finite $G$-set and $\Diam \in \Diam{_G^I}$. For $\Diam$ normal, the Hurwitz stack $\Mgm[g][I][G]_{\Diam}^{rat}$ is obtained as:
				\[\xymatrix{
					\Mgm[g][I][G]_{\Diam}^{rat}\ar[r]\ar[d]&\Hurw[G]^{rat}\ar[d]\\
					\M_g[G]_{\bar\Diam}\ar[r] & \Hurw[G]. 	
				}
				\]
			\end{proposition}	
			
			\medskip
			
			\begin{remark}\label{rem:HurwCaseCyc}
				In the case of $G\simeq \ZZ/n\ZZ$ cyclic, the definition of (normal) \emph{unordered} $\Diam{_G^I}$-Hurwitz data $\bar\Diam$ is identical to the one of Hurwitz data  $\mathbf{k}\in (\ZZ/n\ZZ)^\nu/\mathfrak{S}_\nu$ of \cite{ColMau1} Def.~3.5. The requirement of considering the $\aut(G)$-classes as in ibid. Rem.~3.6 ii) is dealt with later in terms of $\aut_{\lozenge}(G;H)$-quotient, see \S \ref{subsub:QuotientByAut} and \S \ref{sec:IrrHurwLoc}.
			\end{remark}

		\subsubsection{}\label{subsub:GArtih}
			Let us describe briefly the arithmetic of Hurwitz data and in which sense the normality property of a $\Diam_G^I$-Hurwitz data encodes some arithmetic of $G$-covers.

			\medskip

			Let $C\in \Mgm[g][m][G]_{\bar\Diam}(k)$ be a $G$-cover $\pi\colon X\to Y=X/G$ defined over $k$, and let $R_\pi\subset X$ be the ramification divisor that splits over $K\supset k$. The $\text{Gal}(\bar k/k)$-action on $R_\pi$ translates as follow on $\{\Diam\}_{\Diam_G^I}$: it is given (1) on the $\{G_i\}_I$ by conjugacy, and (2) on the $\{\chi_i\}_I$ by change of primitive $N$-root of unity in $\mu_{|G_i|}(K)$.

			\medskip
			
			Let us further denote by $Y'$ the $G$-quotient minus the branch locus of $\pi$. The monodromy representations provides a correspondence between classes of $G$-covers over $Y'$, and $G$-conjugacy classes of morphisms $\psi_{\pi}\colon\pi_1^{et}(Y')\twoheadrightarrow G$ -- see \cite{DD06} \S 2.4.2. Assuming that $Y'$ admits a rational point $x\in Y'$ gives a splitting of the arithmetic-geometric fundamental exact sequence:
			\[
			1\to \pi_1(Y'\times \spec \bar k) \to\pi_1(Y') \to \textrm{Gal}(\bar{k}/k)\to 1 
			\]
			and $\psi_{\pi}$ gives a morphism $\phi\colon\text{Gal}(\bar k/k)\twoheadrightarrow G$: under which $\text{Gal}(\bar k/k)$ acts on the class of $G$-covers. Under this $\text{Gal}(\bar k/k)$-action, the ramification divisor $R_{\pi}$ of $\pi$ is \emph{globally} invariant and $\phi$ reads as the $\text{Gal}(\bar k/k)$-action on the generic unramified fiber over $x$ of the cover -- see Ibid \S 2.9.
			
			\medskip
			
			\pagebreak[4]
			
			In terms of moduli spaces $\Mgm[g][m][G]^{rat}_{\Diam}$ associated to a \emph{normal} Hurwitz data $\Diam$, one recovers via $\phi$ the previously described $\text{Gal}(\bar k/k)$-action, here on $\{\Mgm[g][I][G]^{rat}_{\gamma.\Diam},\ \gamma\in G\}$ via the $G$-action on $\Diam$. This action satisfies the following property.
			
			\begin{proposition}\label{prop:NormArith}
				For $\Diam$ a normal $\Diam_G^I$-Hurwitz data, the moduli stack $\Mgm[g][I][G]^{rat}_{\Diam}$ is $\Gq$-invariant in $\M_{g,[m]}$.
			\end{proposition}

			This result is immediate following the properties of normal Hurwtiz datas developed in the previous sections. By controlling the $\Gq$-action on the characters $\{\chi_i\}_I$, the normality condition thus provides an a priori finer arithmetic invariant than the global $\Gq$-invariance of $R_{\pi}$, i.e the control of the $G$-isotropy groups of $R_{\pi}$ only.

	\subsection{Functorial Properties and Moduli}
		We briefly present some properties of the Hurwitz stack with rational ramification $\Hurw^{rat}$ and of the Hurwitz data $\Diam$ under change of group and under quotient by automorphisms. In what follows we denote by  $H\triangleleft G$ a normal subgroup of $G$ and by $P=G/H$ the associated quotient.

		\subsubsection{}\label{subsub:ResCores}
			Let $(E,\iota_G)$ be a curve with $G$-action, and let us denote $C=E/H$ and $D=E/G$ the various quotients given by the situation above, which induces a restricted $H$-action $(E,\iota_H)$ on $E$ as well as a corestricted $P$-action $(C,\iota_P)$. 

			\medskip
			
			For a given Hurwitz data $\Diam=\{(G_i,\chi_i)\}_{i \in I}$, one obtains some \emph{restricted $\res(\Diam)=\{(H_j,\chi_j')\}_{j \in J}$ and corestricted $\cores(\Diam)=\{(P_k,\tilde{\chi}_k)\}_{k \in K}$ Hurwitz datas} given respectively by
			\begin{equation*}
				\res(\Diam_G)=\begin{cases}
					H_j=H\cap G_j \\
					\chi_j'=\chi_{j|H\cap G_j}					
				\end{cases} %
				\cores(\Diam_G)=\begin{cases}
					P_k=G_k/G_k\cap H\\
					\tilde{\chi}_k=\chi_k^{|G_k\cap H|}
				\end{cases},
			\end{equation*}
			whose definition is functorial and is similar to those of \cite{BERO07} \S 2.2.2.

			Note that  the indices $J$ defining $\res(\Diam)$ (resp. $K$ defining $\cores(\Diam)$) are taken over the indices of $I\in \Diam$ with non-trivial induced isotropy groups $H_i$ (resp. non-trivial $P_i$). These Hurwitz datas $\res(\Diam)$ and $\cores(\Diam)$  correspond respectively to those of the restricted and corestricted actions $(C,\iota_H)$ and $(D,\iota_P)$ on the curves.

		\subsubsection{}\label{subsub:ModCores}
			Let $\Diam$ be a $\Diam_G^I$-Hurwitz data and let us denote by $g'$ the genus of the quotient as defined by $\Diam$. We now write $\Hurw[H\trianglelefteq G]$ instead of $\Hurw$ for the Hurwitz stack of $G$-curves to remember that a normal subgroup $H\trianglelefteq G$ is fixed.
			
			\medskip
			
			By forgetting the action of $G$, resp. considering the (tame) quotient by $H$, one obtains some restricted and quotient morphisms
			\[
			\Mg[H\triangleleft G]\longrightarrow \Mg[H],\text{ resp. }\Mg[H\triangleleft G]\to \M_{g'}[P],
			\]
			as well as their rational variants on $\Hurw^{rat}$. One deduces from the study above their versions with Hurwitz data, $\Mg[H\triangleleft G]^{rat}_{\Diam}\longrightarrow \Mg[H]^{rat}_{\res(\Diam)}$ and also
			\begin{equation}\label{eq:PhiRat}
				\tilde{\Phi}_\Diam^{rat} \colon \Mg[H\triangleleft G]^{rat}_{\Diam}\to \M_{g'}[P]^{rat}_{\cores(\Diam)}.
			\end{equation}

			\begin{proposition}\label{QuotientIsOpen}
				The morphism
				\[
				\tilde{\Phi}\colon \mathcal M_g[H \triangleleft G]^{rat} \to \mathcal M_{g'}[P]^{rat}
				\]
				is smooth.
			\end{proposition}
			
			\begin{proof} Because the action of $H$ is tame, so that $H$ is reductive, we see that $\tilde{\Phi}$ induces first a surjection at the level of tangent spaces by Proposition 4.6 of \cite{MaugeaisBordeaux}, then at the level of geometric points as it can be checked functorially. The result then follows from the smoothness of  $\mathcal M_g[H \triangleleft G]^{rat}$ and $\mathcal M_{g'}[P]^{rat}$, which can be proven using deformation theory using once more the reductiveness of $P$ and $G$, and from the characterisation of smooth morphisms between smooth schemes (cf. \cite{EGA4.4}, Théorème 17.11.1).	
			\end{proof}

		\subsubsection{}\label{subsub:QuotientByAut}
			Recall that one has a left $\aut(G)$-action on $\Mg[G]$. The effect of an automorphism $\psi \in \aut(G)$ on a $G$-Hurwitz data $\Diam$ attached to a curve $(C,\iota_G)\in \Mg[G]_{\Diam}^{rat}$ is read on the $\psi$-twisted $G$-cover $(\tilde C,\tilde{\iota}_{G})\in\Mg[G]_{\psi(\lozenge)}$ as follow: 

			\begin{enumerate}
				\item the action of $G$ on $I$ and $C$ is obtained through the action of $G$ twisted by $\psi$,
				\item the ramification divisor -- with the $G$-action given via the sections through $I$ -- is unchanged $\textrm{Ram}((C,\iota_G))=\textrm{Ram}((\tilde C,\tilde{\iota}_{G}))$, although its numbering is, together with the action on $I$,
				\item a $G$-isotropy group $G_i\in \lozenge$ is sent to one of $(\tilde C,\tilde{\iota}_{G})$: $\psi(G_i)=G_j\simeq G_I<G$,
				\item  the characters $\chi_i\colon G_i\to \mu_{|G|}$ are uniformly changed by composition by $\psi$: $\psi.\chi_i=\chi_i\circ \psi_{|G_i}$. 
			\end{enumerate}
			In particular if $\psi$ fixes $\Diam$ then $\psi(G_i)=G_i$ for all $i\in I$.

			\medskip
			
			For a given $\Diam\in\Diam_G^I$-Hurwitz data, we denote by \emph{$\aut_{\lozenge}(G)$ the subgroups of $\aut(G)$ that fixes $\lozenge$}, and one defines the \emph{subgroup $\aut_{G/H,\Diam}(G;H)$ of $\aut_{G/H}(G;H)$ fixing $\Diam$} -- see Appendix \ref{sec:Erratum}. In particular, we can now consider the $\Diam$-fixing variant $\Hurw[H \triangleleft G]^{rat}_{\Diam}/\aut_{G/H,\Diam}(G;H)$ of $\Hurw[H \triangleleft G]^{rat}/\aut_{G/H}(G;H)$ and the \emph{relative quotient morphism}:
			
			\begin{equation}\label{eq:relQuotMorp}
				\Phi_{\Diam}^{rat}\colon \Hurw[H \triangleleft G]^{rat}_{\Diam}/\aut_{G/H,\Diam}(G;H)\to \mathcal M_{g'}[G/H]^{rat}_{\textrm{cores}^G_{G/H}(\Diam)}.
			\end{equation}
			whose irreducibility is established in \S \ref{sec:IrrRel}.

\section{Mixed Cohomology and Torsors Extensions}\label{sec:MixCohom}
	Let $X$ be a $S$-scheme endowed with an action of a finite group $P$. We study a certain set of torsors extensions associated to the $P$-torsor $X$ and to some finite abelian groups $H$ as given in the figure below. 

	\begin{wrapfigure}[8]{r}{3cm}\centering \vspace*{-1em}
		$\xymatrix{
			Z \ar@{-->}[d]_H \ar@/^1pc/[dd]^G\\
			X \ar[d]_P\\
			Y
		}$ 
	\end{wrapfigure}
	This is achieved in term of \emph{mixed étale cohomology} following \cite[\S 2.1]{ClassesChern}. Recall that the mixed étale cohomology groups $\mathrm{H}^i_P(X_\et, \mathcal F)$ are defined as the derived functors of $\mathcal F \mapsto \mathrm{H}^0_\et(X, \mathcal F)^P$ for a sheaf of abelian groups with a $P$-action $\mathcal F\in ShAb(X,P)$, and that they are the abutment of two spectral sequences -- see Ibid. and \cite[Chap. V]{TOH57} -- which relates these extensions to their branch loci and class of group extensions respectively via their étale and group cohomologies. 
	
	In what follows $X$ denotes an irreducible $S$-scheme endowed with a faithful action of a finite abstract group $P$ such that $Y = X/P$ exists as a $S$-scheme, the group $H$ is an abelian finite group, and $Z\to X$ denotes a torsor extension as above, also designed as an $Aut_Y(Z)$-torsor when $G$ is not identified.

	\subsection{From Mixed Cohomology to Group Cohomology}
		Let $Z$ be a scheme endowed with a $G$-action denoted $i_Z$. We endow the set of torsors extensions with the equivalence relation given by $G$-equivariant isomorphisms given by $H$-conjugacy. 

		\begin{definition} With the notations and assumptions above, we denote by $\mathcal{T}ors_{X,P}(H, G)$ the set of classes of $(H;G)$-torsors extensions of $P$-torsors $X$:
			\[
			\mathcal{T}ors_{X,P}(H, G) = 
			\left\{\begin{array}{l}
			\textrm{$H$-torsor } Z' \to X  | Z' \to Y \textrm{ is a $G$-torsor and }\\
			Z'/H \cong X \textrm{ as $P$-torsor }
			\end{array}\right\}/\sim
			\]	
			where two $(H;G)$-torsors extensions $(Z' \to X)$ and $(Z'' \to X)$ are equivalent if and only if there exists a $G$-equivariant $Y$-isomorphism $\phi\colon Z'\simeq Z''$ and $h\in H$ such that $\phi\circ i_{Z'}(\gamma) = i_{Z''}(h\gamma h^{-1})$ for all $\gamma \in G$. 
		\end{definition}

		It follows from \cite{MaugeaisBordeaux} Proposition 2.2 that the condition for a scheme $Z\in \H^1_{\et}(X,H)$ to be an $Aut_Y(Z)$-torsor over $Y$ is to be a $P$-equivariant $H$-torsor over $X$, i.e. $Z\in \H^1_{\et}(X,H)^P$; in this case the group $\aut_Y(Z)$ is an extension $G$ of $P$ by the kernel $H$. As shown below, the mixed cohomology presents a refinement of this situation.

		\subsubsection{}  For $\mathcal F\in AbSh(X_{et},P)$, let us consider the first spectral sequence of mixed cohomology of \cite[\S 2.1]{ClassesChern}
			\begin{equation}
				\label{eq:seq1}
				E_2^{p, q} = \mathrm{H}^p(P, \mathrm{H}^q_\et(X, \mathcal F)) \Rightarrow \H^{p+q}_P(X_\et, \mathcal F).
			\end{equation}
			
			For a finite group $H$ with a given $P$-action, we define the constant $P$-sheaf $H_X$ on $X$, which leads to the long exact sequence
			\begin{equation}\label{eq:ExSeqGrp}
				0 \to \H^1(P, H_X)\to \H^1_P(X_\et, H_X) \to \H^1_\et(X, H_X)^P \to \H^2(P, H_X)
			\end{equation}		
			where we identify $\H^1(P, \H^0_\et(X, H_X))=\H^1(P, H_X)$ since $X$ is connected. 

			\medskip
			
			In particular, one recovers that two mixed $P$-covers of $X$, with difference $\alpha \in \H^1_P(X_\et,H)$, giving rise to the same element in $\H^1_{\et}(X,H)^P$ have the same $H$-conjugacy class of splitting extension of $P$. Moreover as an $Aut_Y(Z)$-torsor defines a group extension $G\simeq Aut_Y(Z)\in \H^2(P,H)$ of abelian kernel $H$, the difference of two $P$-invariant $H$-torsors of $\H^1_P(X_{et},H)$ with same image in $\H^2(P,H)$ is actually a $G$-torsor over $Y$ for $G\in \H^2(P,H)$.

		\subsubsection{} Following the discussion above, the exact sequence \eqref{eq:ExSeqGrp} defines an action of $\H^1_P(X_\et, H)$ on the set of $(H;G)$-torsors extensions as defined via the association of a $G\in \H^2(P,H)$, and this action is well-defined on the classes in $\mathcal{T}ors_{X,P}(H, G)$ via $\H^1(P,H)$.

			\begin{theorem}	\label{th:IdentH1-Tors}
				Let $P$ be a finite group, $H$ a finite abelian group, $Y$ a scheme and $X/Y$ a $P$-torsor. For any extension $G$ of $P$ by $Y$, the mixed cohomology group $\H^1_P(X_\et, H)$ acts simply transitively on $\mathcal{T}ors_{X,P}(H, G)$.
			\end{theorem}
			
			This result is a key ingredient in establishing our main Theorem \ref{th:RelativeIrr} as it allows a fine control of the set of $(H;G)$-torsors.

	\subsection{From Mixed Cohomology to Cohomology of the Quotient Space} 
		Let us consider the second spectral sequence converging to the mixed cohomology groups. With the same notations as above, let us denote by $\pi \colon X \to Y=X/P$ the quotient morphism. For $\mathcal F\in ShAb(X,P)$, one has as a special case of the Leray spectral sequence:
		\begin{equation}\label{eq:seq2}
			E_2^{p, q} = \H^p_\et(Y, \R^q\pi_*^P \mathcal F) \Rightarrow \H^{p+q}_P(X_\et, \mathcal F).
		\end{equation}
		When $\pi$ is étale -- i.e. the action of $P$ is free -- one has $\R^q\pi_*^P \mathcal F=  0$ for $q > 0$ and we get by degeneracy:
		\[\H^{p}_P(X_\et, \mathcal F) = \H^p_\et(Y, \pi_*^P \mathcal F). \]
		
		In what follows, we consider a finite abelian group $H$ endowed with a given action of $P$, and we denote by $H_X$ the associated sheaf on $X$. In this case, the sequence \eqref{eq:seq2} leads to a local-global long exact sequence:
		\begin{equation}\label{eq:loc-glob}
			0 \to \H^1_\et(Y, \pi_*^P H_X) \to \H^1_P(X_\et, H_X) \to \H^0_\et(Y, \R^1\pi_*^P H_X) \to \H^2_\et(Y, \pi_*^P H_X).
		\end{equation}
		We establish a lifting property of the global part $\H^1_\et(Y, \pi_*^P H_X)$ of the mixed cohomology group which is a key element in controlling the deformation of $(H;G)$-covers in the proof of Theorem \ref{th:RelativeIrr}.
		
	\subsubsection{}\label{subsubs:Norm} 
			To compute the global part $\H^1_\et(Y, \pi_*^P H_X)$ we consider the following norm morphism: 
			\begin{equation}
				\label{eq:defNorm}
				N_P \colon H_Y \to \pi_* \pi^* H_Y = \pi_* H_X \to \pi_*^P H_X
			\end{equation}
			defined by $N_P(a) = \sum_{\bar g \in P} \bar g.a$. Let us denote by $\mathcal K$ is kernel and by $\mathcal C$ its cokernel. For $y \in Y$ in the \'etale locus of $\pi$ the norm $N_P$ induces an isomorphism in a neighbourhood of $y$. The sheaves $\mathcal K$ and $\mathcal C$ have therefore their support in the branch locus. 
			
			\medskip
			
			We can actually be more specific: let $x \in X$ be a ramified point of $\pi$ and let $y = \pi(x)\in Y$.  Denoting by $P_x$ the stabiliser of $x$ in $P$, we have 
			\[(\pi_* H_X)_y = \bigoplus_{x' \in P.x} H_{X, x'} \text{ thus } (\pi_*^P H_X)_y \cong H_{X, x}^{P_x},\]
			so we identify $\mathcal K_y$ and $\mathcal C_y$ to the kernel and cokernel of the norm morphism $N_P^{loc}\colon H_{Y, y} \to H_{X, x}^{P_x}$. From this we deduce:		
			
			\begin{proposition} Let $x \in X$ be a ramified point of $\pi$ of stabilizer $P_x$ in $P$ and let $y = \pi(x)$. Then there exists isomorphisms
				\begin{equation*}\label{eq:KrnCoKern}
					\mathcal K_y \cong \H^0(P_x, H_{X, x}) \text{ and } \mathcal C_y \cong \H^2(P_x, H_{X, x})
				\end{equation*}
				where $\mathcal K_y$ and $\mathcal C_y$ are the kernel and cokernel of $N_P^{loc}$ above.
			\end{proposition}

		\subsubsection{} Suppose moreover that $X$ is a curve and the action of $P$ is faithfull. 
			Since $\mathcal K$ and $\mathcal C$ are supported in the branch locus of $\pi$ which is discrete, their higher cohomology groups vanish and we obtain a long exact sequence:
			\begin{multline*}
				0 \to \H^0_\et(Y, \mathcal K) \to \H^0_\et(Y, H_Y) \to \H^0_\et(Y, \pi_*^P H_X) \to \H^0_\et(Y, \mathcal C) \to \\
				\H^1_\et(Y, H_Y) \to \H^1_\et(Y, \pi_*^P H_X) \to 0
			\end{multline*}
			
			From this we deduce the main property of this subsection.
			\begin{proposition}\label{prop:LiftingGlobalElement}
				The norm morphism $N_P$ of Eq. \eqref{eq:defNorm} induces a transitive action of $\H^1_\et(Y, H_Y)$ on $\H^1_\et(Y, \pi_*^P H_X)$.
			\end{proposition}
			
			In particular, this reduces the construction of $P$-equivariant $H$-torsors over $Y$ to the construction of $H\triangleleft G$ torsors over $Y$ -- see Corollary \ref{cor:defGlob} for the application to deformation of the global part of $(H;G)$-covers.

		\subsubsection{}
			We now give a group theoretic description of $\H^1_P(X_{et},H_X)$.

				\medskip
				
				In terms of group theory, the same approach that in \S \ref{subsubs:Norm} shows that the fibers at $y\in Y$ of the sheaf $\R^1 \pi_*^P H_X$ identifies to
				\begin{equation}
					\left(\R^1 \pi_*^P H_X\right)_y \cong \H^1(P_x, H),
				\end{equation}
				where $x\in X$ is a lifting of $y$, via the isomorphism $\left(\R^1 \pi_*^P H_X\right)_y\cong\H^1(P, (\pi_* H_X)_y)$. This implies some refinements on the description of $H^1_P(X_{et},H_X)$ in term of kernel/cokernel, that depend on the first and second group cohomology class extensions of $P$ by $H$.
				
				For example, for $H$ and $P$ coprime order cyclic groups with trivial class in $H^2(P,H)$, one obtains that $H^1_P(X_{et},H_X)\simeq H^1(Y,H)$. The other cases of class of $G\in H^2(P,H)$ lead to similar identification of $H^1_P(X_{et},H_X)$ that we do not reproduce here, since we do not make use of these computations.

\section{Irreducibility of the Relative Quotient Morphism}\label{sec:IrrRel}
	Let $G$ be a finite group, $H\triangleleft G$ a subgroup and $P=G/H$ its quotient, and let $\Diam \in \Diam^I_G$ be a Hurwitz data. The goal of this section is to prove the irreducibility of the quotient morphism  $\Phi_{\lozenge}^{rat}$ of Eq.~\eqref{eq:relQuotMorp} when $\lozenge$ has a $P$-\'etale point, i.e when one of the $e_i$ of $Ram(C,\iota)$ has an image modulo $H$ with trivial stabilizer.
	\begin{theoNoSub}\label{th:RelativeIrr}
		Let $G$ be a finite group, $I$ a finite $G$-set, $H\triangleleft G$ be a cyclic group of prime order and let $\lozenge$ be a $\Diam_G^I$-Hurwitz character. If $\lozenge$ has a $P$-\'etale point, then the morphism
		\[
		\Phi_\Diam^{rat}\colon \mathcal M_g[G]^{rat}_{\Diam}/\aut_{G/H, \Diam}(G;H) \to \mathcal M_{g'}[G/H]^{rat}_{\cores(\Diam)}
		\]
		is geometrically irreducible.
	\end{theoNoSub}
	
	Recall that $\aut_{G/H, \Diam}(G;H)$ denotes the subset of $\aut(G)$ that leaves $G/H$ and $\Diam$ fixed -- see \ref{def:autPfix}.
	
	\medskip
	
	The context, along the lines of \cite{MaugeaisBordeaux}, is given by whose of $(H;G)$-torsors, which allows to go back and forth to $G$-curves, as well as to apply the results of \emph{mixed cohomology} on \'etale and group cohomology of the previous section \S \ref{sec:MixCohom}. The proof follows two steps: first a local deformation of the ramification locus -- see Proposition \ref{prop:DefBrLoci} -- up to a global part, then the deformation of the global part in cohomology --  see Corollary \ref{cor:defGlob}. The local deformation relies on the arithmetic properties of torsors and their compactification to construct explicit $P$-equivariant deformations of $H$-torsors up to a global torsor. The global deformation then relies on equivariant generalisations of ideas developed in \cite{COR87}. 
	
	\medskip
	
	Before going into the deformation process, we first recall some properties of the stack $[R^1f_*H]^{G/H}$ that classifies the $G/H$-invariant $H$-torsors.

	\subsection{Stack of Covers and Stack of Torsors}\label{subsub:R1fG}
		Following \cite{MaugeaisBordeaux}, a stack of $G$-equivariant curves can be described in terms of a certain stack of torsors $[\R^1f_*G]$. Indeed, let $\mathcal C\to \mathcal M_{g',[m]}[G/H]$ be the universal $G/H$-equivariant curve of genus $g'$ endowed with an equivariant divisor $\mathcal B$ of degree $m$. Let us write  $f\colon \mathcal C\setminus|\mathcal B|\to \mathcal M_{g',[m]}[G/H]$.
		Following ibid. \S 3.2, one attaches to this situation an algebraic stack $[R^1f_*H]^{G/H}$ over $\mathcal M_{g',[m]}[G/H]$ that classifies the local $H$-torsors that are $G/H$-invariants. This stack admits a stratification $[R^1f_*H]^{G/H}_{\mathfrak{g}=g}$ given by the genus $g$ of the fibres -- see ibid. \S 4.2 -- and in our case fits within the diagram -- see ibid \S 4.4:
		\begin{equation}\label{eq:fullDiag}
			\xymatrix{
				\mathcal M[H\triangleleft G]_{\lozenge}\ar[r]\ar[d]  &[R^1f_*H]^{G/H}_{\mathfrak{g}=g}\ar[d] \\
				\mathcal M_g[H\triangleleft G]^{rat}_{\lozenge}/\aut_{G/H, \lozenge}(G; H) \ar@{-->}[ru]^{\Psi}\ar[r]^(0.57){\Phi}& \mathcal M_{g'}[G/H]^{rat}_{\cores[G][G/H](\lozenge)} 
			}					
		\end{equation}		
		where $H\triangleleft G$ is here to recall that we implicitly deal with $(H;G)$-covers. In this diagram, $\aut_{G/H, \lozenge}(G; H)$ denotes a certain subgroup of the automorphisms of $G$ preserving $H\triangleleft G$ -- see Appendix \ref{sec:Erratum} for definition and properties -- and $\Psi$ factorizes the upper-left triangle, cf. Eq. \eqref{eq:AppRfGfact}.
		
		\medskip
		
		This implies the identification of the fibres of $\Phi$ to a space of torsors, or more precisely (see \S 5 of \cite{MaugeaisBordeaux}):
		\begin{proposition}\label{prop:RfH} 
			In the situation above, let us consider $\ZZ/p\ZZ\simeq H\triangleleft G$. Then $[R^1f_*H]^{G/H}_{\mathfrak{g}=g}$ is representable by an algebraic stack and $\Psi$ induces an isomorphism onto its image.
		\end{proposition}
		
		This construction allows to identify equivariant spaces of curves to some particular torsors. In terms of the analytification of these spaces, we are able to construct topological paths on $[R^1f_*H]^{G/H}_{\mathfrak{g}=g}$ that automatically lift to $\mathcal M_g[H\triangleleft G]/\aut_{G/H}(G; H)$ thanks to the isomorphism $\Psi$. We emphasize the importance of the quotient by $\aut_\bullet(G;H)$ in order for $\Phi$ to be a local immersion -- see Appendix \ref{sec:Erratum}.

	\subsection{Deformation of Branch Locus}\label{sub:DefLoc} 
		Let $E$ and $E'$ be two $G$-curves with respective rational ramification $R_E$ and $R_{E'}$, and same $\Diam_G^I$-Hurwitz data $\Diam$. In particular, $E $ and $E'$ have same $G$-quotient $D$. Suppose that $E$ and $E'$ have also isomorphic $H$-quotient $C=\Phi_\Diam^{rat}(E) = \Phi_\Diam^{rat}(E')$ in $\mathcal M_{g'}[P]^{rat}_{\cores(\Diam)}$. Assuming that $\Diam$ has a $P$-étale point, our goal is to deform $E$ into an $G$-equivariant curve with quotient $D$, same Hurwitz data $\Diam$ and with an $H$-branch locus equal to that of $E'$.

		\subsubsection{}\label{subsub:TheIJContext}
			We denote by $J_H$ the indices of $\res(\Diam)$. Writing $R_E=(e_i)_{i \in I}$ and $R_{E'}=(e'_i)_{i \in I}$, we consider 
			\[
			Z = E \setminus \cup_{i \in J_H} \{e_i\}\text{ and }Z' = E' \setminus \cup_{i \in J_H} \{e'_i\}
			\] 
			with $X = Z/H$ and $X' = Z'/H$, so that $Z$  (resp. $Z'$) is exactly the étale loci of $\pi_H\colon E \to C$ (resp. $\pi'_{H} \colon E' \to C$). In particular,  $Z \to X$ and $Z' \to X'$ are both $H$-torsors.
			
			\begin{wrapfigure}[9]{r}{4cm}
				\vspace*{-1em}
				\centering
				\xymatrix@C.3em@R3em{
					Z \ar[d]& & Z'\ar[d] & E,\, E'\ar[d]^H\ar@/^1.7pc/[dd]^{\pi^{({\scriptstyle '})}_{G}}\\
					X\ar[dr]\ar@{}[r]|-*[@]{\supset}& \tilde X \ar@{}[r]|-*[@]{\subset}& X'\ar[dl]& C\ar[d]^P_{\pi}\\
					& Y & & D\\
				}
			\end{wrapfigure}
			
			\medskip
			
			Let $I_H^P\subset J_H$ be the $P$-étale indices of the restricted $H$-cover, i.e. $I_H^P = \{i \in J_H, G_i \subset H\}$. By definition $I_H^P$ is exactly the set of indices of $J_H$ for which $e_i\in R_E$ and $e'_i\in R_{E'}$ are sent to étale points in $C \to C/P$. 	Notice that $I_H^P$ is naturally endowed with a \emph{free} $P$-action because if $i \in I_H^P$, then $G_i = H$ as $H$ is of prime order.

			Let $R_{E}^P=\{\pi_H(e_i)\}_{i \in I_H^P} $ and $R_{E'}^P=\{\pi'_{H}(e'_i)\}_{i \in I_H^P}$, which are the ramification points of $C \to C/P$, since $\pi_H(e_i) = \pi'_H(e'_i)$ for $i \in I \setminus I_H^P$. In the deformation process, the points of $R_{E}^P $ and $R_{E'}^P$ are thus the only ones that have to be moved, and we assume $I_H^P \neq \emptyset$ accordingly. This deformation of $R_{E}^P $ to $R_{E'}^P$ is achieved using torsors and their compactification.

			\medskip
			
			Finally, the schemes $X$ and $X'$ can actually be embedded into $C$ in a $P$-equivariant way since $E/H = E'/H = C$. We then consider $\tilde X = X \cap X'$, which inherits a $P$-action, and identifies $Z$ and $Z'$ are elements of the mixed cohomology group $\H^1_P(\tilde X, H)$ by Theorem \ref{th:IdentH1-Tors}.

		\subsubsection{}\label{subsub:CompGTors}
			We now prove the existence of a deformation of $E$ whose difference with $E'$ is actually in $\H^1(Y, \pi_*^P H)$, i.e. whose image in $\H^0(Y, \R^1\pi_*^P H)$ is zero according to exact sequence \ref{eq:loc-glob}. This result is obtained by deforming the branch locus in $C$ and relies on the existence of compactification of torsors in family as in \cite{MaugeaisBordeaux}. 
			
			\begin{proposition}\label{prop:DefBrLoci}
				Let $E$ and $E'$ be two $G$-curves, with $H\triangleleft G$, and under the assumptions above, in particular with a $P$-étale point. Then there exists an algebraic deformation $\hat E$ of $E$ such that the difference of the $H$-torsors induced by $E'$ and $\hat E$ belongs to in $\H^1(\tilde X,\pi_*^PH)$.
			\end{proposition}
			
			In other words, the difference of $E'$ and $\hat E$ is global. We denote by $\mathcal E/\mathcal S$ this algebraic deformation, whose construction occupies the rest of this section and is done in two steps: first with respect to the Hurwitz data, and then to the class of $G$.
			
			\medskip

			We can suppose that there exists a $\hat{\imath} \in R^P_E\setminus R^P_{E'}$. Otherwise, $R^P_E= R^P_{E'}$ implies that $X = X'$, and the equality of the Hurwitz data implies that the difference is already in $\H^1(\tilde X,\pi_*^PH)$ thanks to the exact sequence \eqref{eq:loc-glob}. Let us move $t_0=\hat{\imath}\in R^P_E\subset C$ to a point of $R^P_{E'}$ by keeping fixed the points that are already in $R^P_E\cap R^P_{E'}$. 
			
			To do so, we first build a family $\hat Z \to C \times \tilde D^\circ$ of $H$-torsors that:
			\begin{enumerate}
				\item is $P$-invariant;
				\item is branched over each $e_i\times D^\circ$ with ramification data $(G_i, \chi_i)$, for $i \in I_H \setminus (G.\hat{\imath})$;
				\item is branched over $P.\pi_H(e_{\hat{\imath}})$ over one fibre and over $P.\pi'_H(e'_{\hat{\imath}})$ over another;
				\item is branched over a constant divisor denoted $\hat\infty$.
			\end{enumerate}
			
			For this, let us fix a point $\infty \in D$ distinct from all the $\pi_G(e_i)$ and $\pi'_G(e'_i)$, and let us define
			\[
			\hat{D}=D\times D^\circ\text{, where we set } D^\circ = D \setminus \bigg( \big\{ \pi_G(e_i), i \in I_H\setminus G.\hat\imath\big\} \cup\{\infty\} \bigg),
			\] 
			which is naturally a family of curves over $D^\circ$. We then consider the relative Cartier divisor $\hat B$ defined through the diagonal $D^\circ \to D \times D^\circ$, which by definition does not meet the pulback $\hat\infty$ of $\infty$ along $\hat D \to D^\circ$. 
			
			It follows from Kummer theory -- see \cite{ColMau1} \S 4.1 -- that, up to a an étale base change $\tilde D^\circ \to D^\circ$, there exists a family of $H$-torsors $\hat Y_H$ over $(\hat D \setminus (\hat B \cup \hat\infty)) \times_{D^\circ} \tilde D^\circ $ that has Hurwitz data equal to $(G_{\hat{\imath}},\chi_{\hat{\imath}})$ over $B$ and equal to $(G_{\hat{\imath}}, \chi_{\hat{\imath}}^{-1})$ over $\infty\times D^\circ$. Up to another base change, let us fix some liftings $f_{\hat\imath}$ (resp. $f'_{\hat\imath}$) of $\pi_G(e_{\hat\imath})$ (resp. of $\pi'_G(e'_{\hat\imath})$) in $\tilde D^\circ$. By pulling back along $\hat C = C\times \tilde D^\circ \to D\times \tilde D^\circ$, we get an \emph{equivariant family $\hat X_H$ of $H$-torsors which are $P$-equivariant, and that are ramified along the pullbacks $\hat B_X$ and $\hat\infty_X$ of $\hat B$ and $\hat \infty$}.
			
			Notice that the fibre $(\hat B_X)_{f_{\hat\imath}}$ is naturally identified to $\sum_{i \in G.\hat\imath} e_i$, and that $(\hat B_X)_{f_{\hat\imath}}$ is identified to $\sum_{i \in G.\hat\imath} e'_i$; the branch points are thus the right ones, and so are the ramification data by construction. The map $\hat X_H \to \hat D$ is however $P \times H$-equivariant only, i.e. with respect to the trivial extension of $P$ by $H$ \emph{only}, and one must still recover the proper class of $G$. We refer to Figure \ref{fig:DefDiag} for a summary of the construction at this stage.
			
			\begin{figure}
				$\xymatrix{
					\hat X_H \ar[rd] \ar[d]\ar@{}[rdd]|-{\square} \\
					\hat Y_H \ar[rd] \ar[d] &\hat C \ar[d] \\
					\hat D \setminus(\hat B \cup \hat\infty)  \ar@{^(->}[r] & D \times D^\circ = \hat D \ar[d]\\
					& D^\circ}$
				\caption{Deformation of $(H;G)$-torsor}
				\label{fig:DefDiag}
			\end{figure}                         
			
			\medskip
			
			The next step is thus to go from $H \times P$ to $G$. For this, we define a family $\hat Z$ over $\hat X$ of $H$-torsors, where 
			\begin{equation}
				\hat Z = \left(Z-(\hat X_H)_{f_{\hat\imath}}\right) \times \tilde D^\circ + \hat X_H \text{ and } \hat X = \hat C \setminus \left(\hat B_X \cup \hat \infty_X \cup \left(\cup_{i \in I_H} e_i \times \tilde D^\circ\right)\right).
			\end{equation}
			This family is obtained by base change and addition of $H$-torsors over $\hat X$, and it satisfies the previous conditions $(1)-(4)$. By Theorem \ref{th:IdentH1-Tors}, this $P$-invariant $H$-torsor defines an element $H^1_P(\tilde{X}_{et},H)$, and thus inherits an action of a group $\hat G$ which is an extension of $P$ by $H$. Since the restriction of this torsor over $f_{\hat\imath}$ equal to $Z$, we obtains $\hat G = G$.
			
			\medskip

			The $\mathcal S$-family $\mathcal E$ of deformation is finally obtained as \emph{a $G$-equivariant compactification of $G$-torsor}. Since the Hurwitz data are constant, so is the genus of the compactification of the fibres of $\hat Z \to \tilde D^\circ$. Up to an alteration $\mathcal S \to \tilde D^\circ$, the existence of such a $G$-equivariant compactification $\mathcal E \to \mathcal S$ of $\hat Z \times_{\tilde D^\circ} \mathcal Z$ follows from \cite{MaugeaisBordeaux} Theorem 4.7. 
			
			Regarding the Hurwitz datas, there exists moreover liftings $\tilde f_{\hat\imath}$ and $\tilde f_{\hat\imath}'$ of  $f_{\hat\imath}$, resp. $f_{\hat\imath}'$, in $\mathcal S$, such that $\mathcal E_{\tilde f_{\hat\imath}} = E$. The Hurwtiz data of $\mathcal E_{\tilde f_{\hat\imath}'}$ over $\pi_H(e_{\hat\imath})$ and of $E'$ at $\pi'_H(e'_{\hat\imath})$ are thus equal, while the Hurwitz datas at the other poins of $R^P_E$ are not changed. Repeating this process over all the points over which the Hurwitz data of $E$ and $E'$ differ, one builds a similar $\mathcal S$-deformation $\mathcal E$ of $E$ into a curve $\tilde E$. 
			
			\medskip
			
			The difference between the $H$-torsors induced by $\tilde E$ and $E'$ is thus global in $\H^1(\tilde X,\pi_*^PH)$ and this concludes the proof of Proposition \ref{prop:DefBrLoci}.

	\subsection{Deformation of the Global Part}\label{sub:DefGlob}
		In the situation of $G$-covers with a $P$-étale point -- see \S \ref{subsub:TheIJContext} for the context --, we establish a deformation result of the global part of the $(H;G)$-torsors. Joined to the local deformation result of Proposition \ref{prop:DefBrLoci}, this concludes the proof of Theorem \ref{th:RelativeIrr}.

		\medskip
		
		The following is a direct application of \S \ref{subsub:CompGTors}, and relies on the analytification of the spaces. 
		
		\begin{proposition}\label{cor:defGlob}
			Let $C$ be a complete smooth curve, $X \subset C$ be an open subscheme. Let $Z$ and $Z'$ be two $(H;G)$-torsors over $X$ with same local datas such that $Z/H \cong Z'/H$ as $P$-torsors. If there exists $x\in C$ with trivial stabiliser in $P$ and over which the completion $E$ of $Z$ and $E'$ of $Z'$ are branched, then $E$ can be deformed continuously and $G$-equivariantly to $E'$ in $\mathcal M_g[H\triangleleft G]/\aut_{G/H}(G; H)$.
		\end{proposition}
		
		Under this assumptions, notice that $E$ and $E'$ have in particular a $P$-\'etale point.

		\begin{proof}
			According to the exact sequence \eqref{eq:loc-glob}, the difference $c = Z-Z' \in H^1_P(X,H_X)$ being global means that $c\in H^1_{\et}(Y,\pi_*^PH_X)$, and it follows from Proposition \ref{prop:LiftingGlobalElement} that $c$ comes from an element $\bar c\in H^1_{\et}(Y,H_X)$. 
			
			\medskip
			
			Let us denote by $\pi_P\colon C\to D$ the $P$-quotient morphism and let us write $y = \pi_P(x)$ that we take as base point of $\pi_1^{top}(D(\mathbb C)^{an},y)$. The Betti-étale comparison isomorphism $\H^1_\et(D, H) \simeq \H^1(D(\mathbb{C})^{an}, H)$ -- see \cite{Milne}, Theorem III.3.12 -- being induced by the morphisms
			\[
			\pi_1^{top}(D(\mathbb{C})^{an}, y) \twoheadrightarrow \H^1(D(\mathbb{C})^{an}, \mathbb{Z}) \twoheadrightarrow \H^1_\et(D, H),
			\]
			the $H$-torsor $\bar c$ over $Y$ lifts to a topological loop $\gamma\in \pi_1^{top}(D(\mathbb{C})^{an}, y_{i})$.
			
			Using the notations of \S\ref{subsub:R1fG}, the loop $\gamma$ defines a path of torsors in $\gamma'$ in $[R^1 f_* H]^{G/H}_{\mathfrak g = g}$ by considering  a family of $G/H$-equivariants $H$-torsors that is branched over the (moving) path $\gamma\subset C$ and over all the other (fixed) branched points coming from $E\to D$ -- they  are equal to those of  $E'\to D$ as the difference between $Z$ and $Z'$ is global. Finally, the path $\gamma'$ can be lifted into  $\mathcal M_g[H\triangleleft G]/\aut_{G/H}(G; H)$ by Proposition \ref{prop:RfH} and it links $E$ to $E'$.
			
		\end{proof}
		
				The proof of Theorem \ref{th:RelativeIrr} is now straightforward under the assumptions of $H$ to be cyclic of prime order and of $\Diam$ to have a $P$-étale point: for a given $E\in \mathcal M_g[G]^{rat}_{\Diam}/\aut_{P, \Diam}(G;H)$, Proposition \ref{prop:DefBrLoci} gives an algebraic deformation of the ramification locus in $C=E/H$ up to a global torsors, which in turn admits a topological deformation by Corollary \ref{cor:defGlob}, hence $\Phi^{rat}_{\Diam}$ is geometrically irreducible.
		
				\bigskip
		
			\begin{remark}
				In the case of Hurwitz data without $P$-étale point, the fibres of $\Phi_\lozenge^{rat}$ are discrete thus not geometrically irreducible; the fact that $H$ does not contain any $G$-isotropy group does not allow the deformation of the global part of the $G$-cover. The morphism $	\tilde{\Phi}_\Diam^{rat} \colon \Mg[H\triangleleft G]^{rat}_{\Diam}\to \M_{g'}[P]^{rat}_{\cores(\Diam)}$ is nevertheless still flat by Proposition~\ref{QuotientIsOpen}.
			\end{remark}

\section{Irreducibility of rational Special Loci}\label{sec:IrrHurwLoc}
	We establish the main result of this paper, which since the normalisation of $\M_{g}(G)$ identifies to $\M_{g}[G]/\aut(G)$, is a special loci version of Theorem~\ref{th:RelativeIrr} for the \emph{rational special loci} $\M_g(G)^{rat}_\Diam=\M_{g}[G]^{rat}_\lozenge/\aut_\lozenge(G)$. To this intent, we introduce the following terminology: 
	
	\begin{defiNoSub} Let $G$ be a finite group, $H\triangleleft G$ and $\lozenge$ be a $G$-Hurwitz data.
	\begin{enumerate}
		\item the group $G$ is said to be \emph{$\aut_{\lozenge}$-liftable with respect to $H\triangleleft G$ and $\lozenge$} (resp. \emph{aut-liftable})  if the morphism $\aut_\Diam(G;H)\to \aut_{\text{cores}_{G/H}^G(\Diam)}(G/H)$ (resp. $\aut(G;H)\to \aut(G/H)$) is surjective;
		\item  the Hurwitz data $\Diam$, is said to be \emph{$H$-irreducible} if the corresponding quotient stack $\M_{g'}[G/H]_{\text{cores}_{G/H}^G(\lozenge)}^{rat}/\aut_{\text{cores}^G_{G/H}(\lozenge)}(G/H)$ is irreducible.
	\end{enumerate}
	\end{defiNoSub}

	The main result of this paper is now:
	\begin{theoNoSub}\label{th:IrredGeneral}
		Let $G$ be a finite group and $\Diam$ be a $G$-Hurwitz data. Assume that $G$ is $\aut_{\lozenge}$-liftable with respect to some $H\simeq \ZZ/p\ZZ\triangleleft G$ such that $\Diam$ has a $G/H$-étale point. If $\Diam$ is moreover $H$-irreducible, then $\M_g(G)^{rat}_\Diam$ is irreducible.
	\end{theoNoSub}

	After some discussion and illustration of the aut-liftability properties, we give the proof of the main theorem, then provide an application to the irreducibility of various {rational special loci} $\M_g(G)^{rat}_\Diam$ under a certain $\aut_\lozenge$-solvable heuristic -- see \S \ref{subsub:AutLozHeur} and Remark~.\ref{rem:MoreIrrLoci}.
	
	\subsection{Lifting Automorphisms with Hurwitz Data} We review various examples, counter-example and properties to the liftability of morphisms to automorphisms with or without Hurwitz data conditions.
	
		\newcommand{\Hom}{\textrm{Hom}}
			
		\subsubsection{Aut-liftablity: a Panorama} A straightforward argument provides two canonical classes of aut-liftable groups.
				
		\begin{proposition}\label{Prop:autLiftDnAb}\mbox{}
			\begin{enumerate}
				\item Finite cyclic groups, and groups of the form $(\ZZ/p\ZZ)^N$ are aut-liftable with respect to their factors;
				\item Dihedral groups $D_{2n}=\ZZ/n\ZZ\rtimes \ZZ/2\ZZ$ are aut-liftable with respect to their rotation groups $H\simeq \ZZ/n\ZZ$.		
			\end{enumerate}
		\end{proposition}
		
	\medskip
	
	We provide 2 examples which illustrate the importance of the choice of $H\triangleleft G$ for aut-liftability in the case of abelian and split metacyclic groups.

	\begin{example}\label{Ex:Cex}\mbox{}
		\begin{enumerate}
			\item \label{Ex:AbNotLift}
			Let $G=\ZZ/p\ZZ\times \ZZ/p^2\ZZ$ and $H=\langle (0,\ p)\rangle$. 
			
			\noindent Then the automorphism $\bar\varphi\in Aut(P)$ defined by the matrix $\begin{pmatrix}0 & 1\\ 1 & 0 \end{pmatrix}$ is not liftable to $G$: otherwise, such a lifting $A=(a_{ij})$ leads to $pA(1, 0)=(0, 0)$ as $(1,0)$ is of order $p$ in $G$, while $pa_{12}=p\neq 0$ in $\ZZ/p^2\ZZ$.
			
			\item \label{Ex:MCSNotLift}	
			Let $G= \ZZ/7\ZZ \rtimes_\psi \ZZ/3\ZZ$ defined by the representation $\psi\colon \ZZ/3\ZZ \to Aut(\ZZ/7\ZZ)\simeq \ZZ/6\ZZ$ given by $\psi(1)=2$, and let us consider $H=\ZZ/7\ZZ$. 
			
			\noindent It then follows from \cite{Curran} Theo.~1 that a general automorphism $\phi\in \aut(G;H)$ is of the form
			\[
			\phi= \begin{pmatrix} \alpha & \beta \\ 0 & \delta\end{pmatrix} 
			\text{ with }
			\begin{cases}
			\alpha \in \aut(\ZZ/ 7\ZZ), & \text{satisfying: }\\
			\beta \colon \ZZ/3\ZZ \to \ZZ/7\ZZ, & \beta(kk')=\beta(k) \beta(k')^{\varphi(\delta(k))}\\ 
			\delta \in \aut(\ZZ/3\ZZ)& \alpha(h^{\psi(k)}) = \alpha(h)^{\psi(\beta(k)\delta(k))}.
			\end{cases}
			\]
			
			where $\delta$ is a given automorphism of $P=\ZZ/3\ZZ$ that we want to lift to $G$.
			
			Since $\beta$ is by definition trivial, the only remaining condition is
			\[\alpha(h^{2k}) = \alpha(h)^{2\delta(k)}\]
			but since $\alpha$ is an automorphism, we get that $2k = 2\delta(k)$ for all $k$, that is $\delta = Id$.
		\end{enumerate}
	\end{example}			

	 In particular Example \ref{Ex:Cex} \ref{Ex:MCSNotLift} shows even automorphisms that are liftable as homomorphism are not always trivially aut-liftable.

	 \begin{remark}\label{rem:GrpCohomIsNoGood}
	 	A group cohomology argument for central group extension and the Schur-Zassenhaus establish the aut-liftablity of non-abelian groups that are product of $2$ coprime factors, one being cyclic. This property will not be used in the rest of the paper.
	 \end{remark}
 
	 \subsubsection{Aut-liftability vs $\aut_{\Diam}$-liftablity}
 	
	 	The characterization of the $\aut(G)$-action on Hurwitz data given in \S \ref{subsub:QuotientByAut} and the injectivity of the characters provide a criterion for an automorphism of $G$ to fix a $\lozenge$ in terms of the isotropy subgroup. This is used in the following proposition to provide a $Aut_\lozenge$-liftable criterion in terms of $R=\langle G_i\rangle_I<G$.
	 	
	 	\begin{proposition}\label{Prop:liftRamDiam} Let $G$, $H\triangleleft G$, $P=G/H$, $\lozenge=\{G_i,\chi_i\}_I$ and $R<G$ be as above, and assume that $\bar\varphi\in\aut_{\cores[G][P](\Diam)}(P)$ lifts to $\varphi\in\aut(G;H)$. Then $\varphi$ fixes $\Diam$ if and only if $\varphi_{|R}=id$. In particular, 
	 		\begin{enumerate}
	 			\item if the image of $R$ generates $P$, then $\bar\varphi$ lifts trivially to $\varphi=Id\in\aut_{\Diam}(G)$. 		
	 			\item if $G$ is abelian, then $\bar\varphi$ lifts canonically to a $\varphi\in\aut_{\Diam}(G)$. 
	 		\end{enumerate}
		\end{proposition}

		\begin{proof}
			The proof is straightforward by noting that the injectivity of the $\chi_i$ forces $\varphi_{|G_i}=Id$, then by writing $G=H\oplus \tilde R\oplus \tilde G$ with $R=H\oplus \tilde R$.
		\end{proof}
	 
	 	To control the a priori non-trivial action of $\aut(G;H)$ on the $G$-Hurwitz data is one of the motivation for working under the $\aut_\lozenge$-solvable Heuristic of \S \ref{subsub:AutLozHeur}, see also below.

	\subsection{Proof of the Rational Irreducibility} 
		Under our $\aut_\lozenge$-liftability and $H$-irreducibility assumptions, the main theorem follows from the relative irreduciblity result of Theorem~\ref{th:RelativeIrr}. Because $G$ is by assumption $\aut_\lozenge$-liftable with respect to $H$ and $\lozenge$, one first obtains a Cartesian diagram:
		\begin{equation}\label{diag:QuotientInduction}
		\begin{aligned}
			\xymatrix{
				\mathcal M_g[H \triangleleft G]^{rat}_{\lozenge} / \aut_{P,\lozenge}(G; H) \ar[r]^(0.6){\Phi_\Diam^{rat}}\ar[d] & \mathcal M_{g'}[P]_{\textrm{cores}^G_P(\lozenge)} \ar[d]\\
				\mathcal M_g[H\triangleleft G]_{\lozenge}^{rat} / \aut_{\lozenge}(G; H)  \ar[r] & \mathcal M_{g'}[P]_{\textrm{cores}^G_P(\lozenge)}^{rat}/\aut(P)}
			\end{aligned}
		\end{equation}
		where one has identifies $\left(\mathcal M_{g'}[H \triangleleft G]^{rat}_\lozenge/\aut_{P, \lozenge}(G;H) \right)/\aut_{\lozenge}(G; H)$ with the quotient in the left bottom corner. 
		
		By $H$-irreducibility hypothesis the right bottom corner is connected, and as $\Phi_\Diam^{rat}$ is geometrically irreducible by Theorem~\ref{th:RelativeIrr} with a surjective left arrow in \eqref{diag:QuotientInduction}, the bottom morphism is geometrically irreducible by base change property as in \cite{EGA4.2} Proposition 4.5.6 (ii). It thus follows from ibid. Proposition 4.5.13, that $\mathcal M_{g'}[H \triangleleft G]^{rat}_\Diam/\aut_{\Diam}(G; H)$ is geometrically connected, hence irreducible since normal. 
		
		Moreover, since the morphism 
		\[ 
		\mathcal M_{g}[H \triangleleft G]^{rat}_\lozenge/\aut_{\lozenge}(G; H) \to \mathcal M_{g}[G]^{rat}_\Diam/\aut_\Diam(G)
		\]
		is surjective, the stack $\M_g(G)^{rat}_\Diam=\mathcal M_{g}[G]^{rat}_\Diam/\aut_\Diam(G)$ is geometrically irreducible by \cite{EGA4.2} Proposition 4.5.4, which concludes the proof.
	
	\subsection{An $\aut_\lozenge$-solvable Heuristic}\label{subsub:AutLozHeur}
		For $G$ a general finite group and $\Diam$ a general $G$-Hurwitz data, the previous section motivates to work under the following assumption -- see also \S \ref{sub:conclusion}.
		
		\medskip
		
		\noindent\textbf{$\aut_\lozenge$-solvable Heuristic.} {\itshape The group G admits an $\aut_\lozenge$-solvable series with respect to a Hurwitz data $\lozenge$, i.e. a series $G = G_0 \triangleright G_1 \triangleright \cdots \triangleright G_m = \{1\}$ such that for every $k\in\llbracket 1,m\rrbracket$: (i) $H_k=G_{k-1}/G_k\simeq \ZZ/p\ZZ$, and (ii) each $G/G_k$ is $\aut_\lozenge$-liftable with respect to $H_k\triangleleft G/G_k$ -- i.e. $\aut_\lozenge(G/G_k; H_k) \to \aut_{\text{cores}_{G/H_k}^G(\Diam)}(G/G_{k-1})$ is surjective for every $k\in\llbracket 1,m\rrbracket$.}

		\medskip
		
	 	We call such a tower $\{G_k\}_{1\leqslant k \leqslant m}$ an $\aut_\lozenge$-resolution of $G$ with respect to $\lozenge$. Note that (i) implies that this series is in particular of maximal length.

		\medskip
		
	  	Under the $\aut_\lozenge$-solvable Heuristic, let us write $\lozenge_k=\mathrm{cores}^{G/G_{m-k+1}}_{G/G_{m-k}}(\lozenge)$ and assume moreover the following two properties to be satisfied: 
	  	\begin{itemize}
	  		\item[$(H_{et})$] that $\lozenge$ induces a sequence $\{(\lozenge_k,H_k)\}_{1\leqslant k \leqslant m}$ so that each $\lozenge_k$ has a $G/G_{k-1}$-étale point
	  		\item[$(H_{quot})$] that the last step $\lozenge_{m-1}$ is $H_{m-1}$-irreducible.
	  	\end{itemize}
		Examples of the latter includes when $\M(G/G_1)$ is a cyclic special loci, i.e. $G/G_1\simeq\langle \gamma\rangle$ as in \cite{ColMau1}. In this situation, an \emph{iterative application of Theorem~\ref{th:IrredGeneral}}, see also Fig.~\ref{fig:ItIrr}, leads to:
		
		\medskip
			
		{\centering
		\begin{minipage}{.9\textwidth}\itshape
		Under the $\aut_\lozenge$-solvable heuristic for $G$ with respect to a $G$-Hurwitz data $\lozenge$, if $\lozenge$ satisfies the properties $(H_{et})$ and $(H_{quot})$, then the rational loci $\M_g(G)^{rat}_\Diam$ is irreducible.
		\end{minipage}  
		\par}
	
		\bigskip

		On the whole, the $(H_{et})$ conditions for general $G$-Hurwitz data, and the $(H_{quot})$ irreducibility condition are both highly impracticable to check for general groups, and thus require some case-by-case and explicit check -- see types $EA$, $C$ and $D^{rot}$ of \S \ref{subsub:CorAbDi} for some examples.
		
		\medskip
		
		\pagebreak[4]
		
		{\addtocounter{figure}{-1}
			\captionof{figure}{Irreducibility of $\M(G):=\M_g(G)^{rat}_\Diam$ under the $\aut_\lozenge$-solvable heuristic and the $(H_{et})$-$(H_{quot})$ hypothesis for $G$ and $\lozenge$.} \label{fig:ItIrr}
			\begin{equation*}\footnotesize
			\begin{tikzcd}[row sep=4pt,column sep=5pt,every matrix/.append style={name=m},execute at end picture={
				\begin{scope}[on background layer]
				\node[rectangle,rounded corners,fill=black!10,minimum width=1cm, fit=(m-1-1) (m-7-1)](c1) {};
				\end{scope};
				\begin{scope}[on background layer]
				\node[rectangle,rounded corners,fill=black!10,minimum width=1cm, fit=(m-1-3) (m-7-3)](c1) {};
				\end{scope};
				\begin{scope}[on background layer]
				\node[rectangle,rounded corners,fill=black!10,minimum width=1cm, fit=(m-1-5) (m-7-5)](c1) {};
				\end{scope};		
			}]
			G \arrow[r,twoheadrightarrow]& G/G_{m-1} \arrow[r,twoheadrightarrow]& G/G_{m-2} \arrow[r,twoheadrightarrow] & ...\arrow[r,twoheadrightarrow] & G/G_2 \arrow[r,twoheadrightarrow] & G/G_{1}\\
			G_1 \arrow[draw=none]{u}[sloped,auto=false]{\leqslant}\arrow[r,twoheadrightarrow]& G_1/G_{m-1} \arrow[r,twoheadrightarrow]\arrow[draw=none]{u}[sloped,auto=false]{\leqslant}& G_1/G_{m-2} \arrow[r,twoheadrightarrow]\arrow[draw=none]{u}[sloped,auto=false]{\leqslant} & \cdots\vdots\cdots \arrow[draw=none]{u}[sloped,auto=false]{\leqslant} \arrow[r,twoheadrightarrow] &  H_{m-1}=G_1/G_{2} \arrow[draw=none]{u}[sloped,auto=false]{\leqslant}& \{1\} \arrow[draw=none]{u}[sloped,auto=false]{\leqslant}\\
			\vdots\cdots \arrow[draw=none]{u}[sloped,auto=false]{\leqslant}& \cdots\vdots\cdots \arrow[draw=none]{u}[sloped,auto=false]{\leqslant} & \cdots\vdots\cdots \arrow[draw=none]{u}[sloped,auto=false]{\leqslant}& \cdots\vdots\cdots \arrow[draw=none]{u}[sloped,auto=false]{\leqslant}&  \{1\}\arrow[draw=none]{u}[sloped,auto=false]{\leqslant} & \\
			G_{m-2} \arrow[draw=none]{u}[sloped,auto=false]{\leqslant} \arrow[r,twoheadrightarrow]& H_2=G_{m-2}/G_{m-1} \arrow[draw=none]{u}[sloped,auto=false]{\leqslant} & \{1\} \arrow[draw=none]{u}[sloped,auto=false]{\leqslant}&  &  & \\
			H_1=G_{m-1} \arrow[draw=none]{u}[sloped,auto=false]{\leqslant}& \{1\} \arrow[draw=none]{u}[sloped,auto=false]{\leqslant}&  &  &   & \\
			\{1\} \arrow[draw=none]{u}[sloped,auto=false]{\leqslant}& &   &  &   & \\
			\M(H_1\triangleleft G) \arrow[r]& \M(H_2\triangleleft G/G_{m-1}) \arrow[r]& \M(H_3\triangleleft G/G_{m-2}) \arrow[r] & \ldots \arrow[r] & \M(H_{m-1}\triangleleft G/G_2) \arrow[r] & \M(G/G_{1})\\
			\end{tikzcd}
			\end{equation*}
		}

		\medskip
			
		It further appears that the two questions of the \emph{irreducibility of the rational special loci $\M_g(G)^{rat}_\Diam=\M_{g}[G]^{rat}_\Diam/\aut_\Diam(G)$ and of the irreducibility of the normal special locus $\widetilde{\M}_g(G)_{\bar\Diam}\simeq \M_{g}[G]_{\bar\Diam}/\aut_\Diam(G)$} are of two complementary arithmetico-geometric nature: for the former, dihedral groups provide some canonical examples of $\aut_\lozenge$-liftability, while Example \ref{Ex:Cex} \ref{Ex:AbNotLift} shows a non-canonical obstruction for abelian group; for the latter, $\aut_\lozenge$-liftable abelian groups automatically provide some irreducible special loci $\M_g(G)_{\bar\Diam}$, while dihedral group show some general descent obstruction from $\M_g(G)^{rat}_\Diam$ to $\widetilde{\M}_g(G)_{\bar\lozenge}$: \emph{a dihedral loci is either irreducible or has rational descent from $\M_g(D_{2p})^{rat}_\Diam$ to $\M_g(D_{2p})_{\bar\Diam}$ } -- see \S \ref{subsub:CorAbDi} and Remark~\ref{rem:D2nIrr}.

		\begin{remark}\label{rem:MoreIrrLoci}
			This approach, under the same hypothesis, implies the irreducibility of each quotient loci $\M(G/G_k)$ for $1\leqslant k\leqslant m-1$. In the case where every $G/G_{m-k}$ is also a \emph{split extension} of $H_{k+1}$ by $G/G_{m-k-1}$, one obtains up to isomorphism a lattice of irreducible subloci $\M(G/G_{m-k-1})<\M(G/G_{m-k})$. 
		\end{remark}	

\section{Cyclic, Elementary Abelian, Dihedral and Irreducible Special Loci}\label{subsub:CorAbDi}
	We provide some elementary blocks of irreducible loci for elementary abelian, cyclic, and dihedral groups that support the $\aut_{\lozenge}$-solvable heuristic. We further illustrate the dual geometric and arithmetic nature of the irreducibility of $\M_g(G)^{rat}_\Diam$ and $\M_{g}(G)_{\bar\Diam}$ in terms of explicit obstruction.

	\subsection{An Irreducible Elementary Abelian Special Loci} 
		Let us consider $G=(\ZZ/p\ZZ)^{N}$, and let $\lozenge$ be any non-empty $G$-Hurwitz data. \emph{The special loci $\M_g(G)_{\bar\lozenge}$ is then irreducible.} 

		\medskip
	
		Indeed, up to a change of factor, one can assume that $H\simeq \ZZ/p\ZZ$ is the first factor of $G$ and that one has $P$-étale ramification. Since $H$ is factor of an abelian group $G$, the latter is \emph{aut-liftable} with respect to $H\simeq \ZZ/p\ZZ$ by Proposition~\ref{Prop:autLiftDnAb} i), while Proposition \ref{Prop:liftRamDiam} ii) gives the $\aut_\lozenge$-liftability of $\Diam$ with respect to $H$.
		
		On the other hand, since $\Diam$ is \emph{$H$-irreducible} for $N=2$ -- since the quotient is a cyclic special loci that is irreducible by \cite{ColMau1} Theorem~4.3 -- one obtains by induction on $N$ that the $\aut_{\Diam}$-solvable heuristic is satisfied, and thus $\M_{g}(G)_\lozenge^{rat}=\M_{g}[G]_\lozenge^{rat}/\aut_\lozenge(G)$ is irreducible.
			
		\medskip
		
		By Proposition~\ref{Prop:normAb}, the Hurwitz data $\lozenge$ is moreover \emph{normal}, and it follows Proposition \ref{prop:OrdUnordNorm} that the normal locus $\widetilde{\M}_g(G)_{\lozenge}\simeq \M_{g}[G]_{\bar\Diam}/\aut_{\Diam}(G)$ is well defined; the special locus $\M_{g}(G)_{\bar\lozenge}$ is thus in turns irreducible.
	
		\medskip
		
		Accordingly, we define:
		
		\medskip
	
			{\centering
			\begin{minipage}{.9\textwidth}\itshape
				\textbf{Type $\bm{EA}$.} For $G=(\ZZ/p\ZZ)^{N}$, $H\simeq\ZZ/p\ZZ\triangleleft G$ first factor, and a Hurwitz data $\lozenge$ that has non-trivial ramification, an elementary abelian loci $\M_g(\ZZ/p\ZZ\triangleleft G)^{rat}_\lozenge$ or $\M_g(\ZZ/p\ZZ\triangleleft G)_{\bar\lozenge}$ is said of type $EA$.
			\end{minipage}  
			\par}
		
		\medskip
		
		We recall that a cyclic special loci $\M_g(\ZZ/n\ZZ)_\Diam^{rat}$ or $\M_g(\ZZ/n\ZZ)_{\bar\Diam}$ -- said of {\bfseries Type $\bm C$} -- is directly irreducible by \cite{ColMau1} and Remark~\ref{rem:HurwCaseCyc}.

	\subsection{Dihedral Special Loci, Irreducibility and Obstruction}\label{sub:TypeDrot} 
		We illustrate how dihedral Hurwitz data provide some examples (1) of descent obstruction from $\M_{g}(D_{2p})^{rat}_\Diam$ to $\M_{g}(D_{2n})_\Diam$ and (2) of irreducible rational special loci $\M_{g}(D_{2p})^{rat}_\Diam$. 	
	
		\subsubsection{} 
			Let $D_{2p}=\ZZ/p\ZZ\rtimes\ZZ/2\ZZ$ with $p\geqslant 5$ with presentation $\langle x,y|x^p=y^2=(xy)^2=1\rangle$, let $\lozenge$ be a $G$-Hurwitz data, and	let $\{G_i\}_I$ denote the $G$-isotropy groups of $\lozenge$. We refer to Proposition~\ref{Prop:normAb}: If one $G_{i}=\langle x\rangle$ is a rotation, then $\lozenge$ is not normal since $G_{i}$ has $(p-1)/2>1$ conjugacy classes of length $2$ and the normality condition $G_{i}\trianglelefteq G$ for $\Diam$ of Eq.~\eqref{eq:HurwCharRel} is thus never fulfilled. On the other hand, if $G_i=\langle xyx^{-1} \rangle\simeq \ZZ/2\ZZ$ is a reflection for all $i\in I$, then this condition holds by unicity of the conjugacy class. The normality of $\Diam$ follows from the second condition given by Eq.~\eqref{eq:HurwCharRel} on the characters $\{\chi_i\colon G_i\to \mu_2\}_I$, which is trivially satisfied. 
			
		\subsubsection{} 
			Let us study the irreducibility with respect to the hypothesis of Theorem~\ref{th:IrredGeneral}. Since $H\simeq \ZZ/p\ZZ$ we have (cyclic) $H$-irreducibility for any Hurwitz data $\lozenge$. The $P$-étale assumption requires that one of the isotropy group $G_{i_0}\simeq \ZZ/p\ZZ$ is a rotation, while one has $\aut_\Diam$-liftablity with respect to any $\Diam$ because of the triviality of $\aut(\ZZ/2\ZZ)$.
	
			\medskip
			
			Accordingly, we define:
			
			\medskip
		
			{\centering
				\begin{minipage}{.9\textwidth}\itshape
					\textbf{Type $\bm D^{rot}$.} For $G=D_{2p}$, $p\geqslant 5$ prime, and $H=\ZZ/p\ZZ\triangleleft D_{2p}$, a dihedral loci $\M_g(\ZZ/p\ZZ\triangleleft D_{2p})^{rat}_\lozenge$ whose Hurwitz data $\lozenge$ contains a rotation is said of type $D^{rot}$.
				\end{minipage}  
				\par} 
			
			\medskip
	
			The rational dihedral loci $\M_g(D_{2p})^{rat}_{\Diam}$ is again in particular irreducible by Theorem~\ref{th:IrredGeneral}.
			
			\begin{remark}\label{rem:D2nIrr}\mbox{}
				\begin{enumerate}
					\item The irreducibility assumptions for $\M_g(D_{2p})^{rat}_{\Diam}$ impose a \emph{non-normal Hurwitz data} which blocks the descent to a potentially irreducible $\M_g(D_{2p})_{\bar \Diam}$;
					\item When $\Diam$ has only rotation type isotropy groups, the type $D^{rot}$ is an example of one of the two cases among three for which the \emph{geometric} invariant $\varepsilon$ of \cite{CLP15} is finer than the Nielsen-numerical one (the other being the case of étale dihedral covers), see below Th.~5.1 ibid, case $2$.
			\end{enumerate}	
			\end{remark}

	\subsection{En Guise de Conclusion}\label{sub:conclusion}
		The $3$ previous types provide some elementary block loci for building irreducible loci $\M_g(G)^{rat}_\lozenge$ and $\M_g(G)$, that should motivates further \emph{explicit studies of Hurwitz data for specific groups}.
	
		\medskip
		
		With the terminology of \S \ref{subsub:AutLozHeur}, let $G$ be a finite group and $\lozenge$ be a $G$-Hurwitz data such that G admits a series $G = G_0 \triangleright G_1 \triangleright \cdots \triangleright G_m = \{1\}$ with induced sequence $\{\lozenge_k,H_k\}_{1\leqslant k\leqslant m}$, $H_k=G_{k-1}/G_k$. If for $1\leqslant k\leqslant m$, the respective loci:
		\begin{enumerate}
			\item $\M_g(H_k\triangleleft G/G_{m-{k-1}})^{rat}$  are either of type $EA$, of type $C$, or of type $D^{rot}$,
			\item $\M_g(H_k\triangleleft G/G_{m-{k-1}})$ are of type $EA$ or of type $C$,
		\end{enumerate} 
		then \emph{the rational loci $\M_g(G)^{rat}_\lozenge$ is irreducible}, resp. the \emph{the special loci $\M_g(G)_{\bar\Diam}$ is irreducible.}
		Indeed, in each case $G$ satisfies the $\aut_\lozenge$-solvable Heurisitic with the $(H_{et})$ and the $(H_{quot})$ properties; in ii) the assumptions ensure descent from $\M_g(G)^{rat}_\lozenge$ to the normalisation $\widetilde{\M}_g(G)_\lozenge$.
		
			\bigskip
		
		{\centering 
			$\star\hphantom{\star}\star \hphantom{\star}\star$
			\par}

\newpage
\begin{appendix}
  	\addtocontents{toc}{\protect\setcounter{tocdepth}{1}}
	\section[Erratum to  ``Quelques calculs d'espaces $R^1f_*G$ sur des courbes'']{Erratum to  ``Quelques calculs d'espaces $R^1f_*G$ sur des courbes'' -- by S. Maugeais}\label{sec:Erratum}

		Let $f\colon \mathcal X\to \mathcal S$ be a morphism of algebraic stacks and $G$ be a smooth group scheme over $S$ and $H\triangleleft G$. In \cite{MaugeaisBordeaux}, we define a certain stack $[R^1f_*\mathcal G]$ over $\mathcal S$ which gives a local classification of $G$-torsors -- denoted $[R^1f_*G]$ in the case of the complement of a relative Cartier divisor in a smooth curve $f\colon \mathcal U\to \mathcal S$ -- and whose $G/H$-equivariant version $[R^1f_*H]^{G/H}$ we show is related to a certain moduli stack $\Mg[H\triangleleft G]$ -- see \S 4.4 \emph{Ibid.} It has to be noted that since the action of $G$ is fixed, \emph{the group $H$ is not a subgroup of the automorphism groups of elements of $\mathcal M_g[H \triangleleft G]$}, and thus one can not consider the  $2$-quotient $\mathcal M_g[H \triangleleft G]\dquot H$ as it is done in \emph{loc. cit.}
	
		Instead of this $2$-quotient by $H$, let us present how the geometric interpretation of $[R^1f_*H]^{G/H}$ is given by forgetting the action of $H$ while keeping that of $G/H$ on the quotient curves.

		\subsection{}
			Let $\aut(G; H)$ denote the subgroup of elements of $\aut(G)$ sending $H$ onto itself. 
			
			\begin{definition}\label{def:autPfix}
				Let $G$ be a finite group and $H\triangleleft G$. We denote by $\aut_{G/H}(G; H)$ the kernel of the  morphism $\aut(G; H) \to \aut(G/H)$.
			\end{definition}
			This group comes with a natural morphism
			\begin{equation}\label{eq:restriction2H}
				\aut_{G/H}(G; H) \to \aut(H)
			\end{equation}
			induced by restriction to $H$ (when viewed as elements of $\aut(G; H)$), which is \textbf{not} injective in general. 
			This groups takes place in an exact sequence
			\begin{equation}
				0 \to \aut_{G/H}(G; H) \to \aut(G; H) \to \aut(G/H)
			\end{equation}
			whose last morphism is not surjective. 
			
			Note that in general, the description of the automorphism group $\aut_{G/H}(G; H)$ is far from being straightforward: see \cite{Wells}. 
			\medskip
			
			Considering the Hurwitz stack with rational ramification $\Hurw^{rat}$ of \S \ref{subsub:RatHurwSt} and the quotient morphism of \S \ref{subsub:ModCores}, one establishes in a similar way to Proposition \ref{QuotientIsOpen}:
		
			\begin{proposition}\label{prop:UnivOpen}
				The morphism
				\[
				\Phi\colon \Hurw[H \triangleleft G]^{rat}/\aut_{G/H}(G;H) \to \mathcal M_{g'}[P]^{rat}
				\]
				is flat, hence universally open.
			\end{proposition}

		\subsection{}
			The $\mathcal M_g[H \triangleleft G]$ stack being endowed with an action of $\aut(G; H)$, one obtains more precisely a factorization of the natural morphism $\Psi$ as below:
			\begin{equation}\label{eq:AppRfGfact}
				\xymatrix{
					\M_{g,g'}[H\triangleleft G]\ar[r]^{\Psi}\ar[rd]& [R^1f_*H]^{G/H}\\
					& \M_{g,g'}[H\triangleleft G]/\aut_{G/H}(G;H)\ar@{-->}[u]
				}
			\end{equation}
			which comes with the corresponding corrections p.~388 and in Proposition~5.2 Ibid. given by replacing the $2$-quotient by $H$ by the quotient by $\aut_{G/H}(G;H)$.
			
			\medskip
			
			The only point that needs verification is actually the point $i)$ of Lemme 4.8 Ibid., because all the others are still valid. This can be rephrased as follows:
			let $\iota_1$ and $\iota_2$ be two actions of $G$ on a curve $C/\spec k$ inducing the same $H$-torsor in a $G/H$-equivariant way. The actions $\iota_1$ and $\iota_2$ are then equal up to an automorphism of $G$ sending $H$ to $H$ and inducing the identity on $G/H$.

\end{appendix}		

\addcontentsline{toc}{section}{References}


\begin{bibdiv}
	\begin{biblist}
		\bib{BERO07}{book}{
			author={Bertin, Jos{\'e}},
			author={Romagny, Matthieu},
			title={{Champs de {H}urwitz}},
			publisher={M{\'e}moire de la SMF},
			date={2011},
			volume={125-126},
			note={\href{https://arxiv.org/abs/math/0701680}{arXiv:math/0701680
					[math.AG]}},
		}
		
		\bib{BRO90a}{article}{
			author={Broughton, S.~Allen},
			title={{Classifying finite group actions on surface of low genus}},
			date={1990},
			journal={Journal of Pure Applied Algebra},
			volume={69},
			pages={233\ndash 270},
			review={\MR{1090743}},
		}
		
		\bib{BROWO07}{article}{
			author={Broughton, S.~Allen},
			author={Wootton, Aaron},
			title={{Finite Abelian subgroups of the mapping class group}},
			date={2007},
			journal={Algebraic {\&} Geometric Topology},
			volume={7},
			pages={1651\ndash 1697},
			review={\MR{2366175}},
		}
		
		\bib{CAT12}{incollection}{
			author={Catanese, F.},
			title={{Irreducibility of the space of cyclic covers of algebraic curves
					of fixed numerical type and the irreducible components of
					{$Sing(\overline{M}_g)$}}},
			date={2012},
			booktitle={{Advances in geometric analysis}},
			series={{Adv. Lect. Math. (ALM)}},
			volume={21},
			publisher={Int. Press, Somerville, MA},
			pages={281\ndash 306},
			review={\MR{3077261}},
		}
		
		\bib{CLP15}{article}{
			author={Catanese, Fabrizio},
			author={L{\"o}nne, Michael},
			author={Perroni, Fabio},
			title={{The irreducible components of the moduli space of dihedral
					covers of algebraic curves}},
			date={2015},
			ISSN={1661-7207},
			journal={Groups Geom. Dyn.},
			volume={9},
			number={4},
			pages={1185\ndash 1229},
			url={http://dx.doi.org/10.4171/GGD/338},
			review={\MR{3428412}},
		}
		
		\bib{ColMau2}{article}{
			author={Collas, Benjamin},
			author={Maugeais, Sylvain},
			title={{On Galois action on stack inertia of moduli spaces of curves }},
			date={2014},
			eprint={1412.4644 [math.AG]},
		}
		
		\bib{ColMau1}{article}{
			author={Collas, Benjamin},
			author={Maugeais, Sylvain},
			title={{Composantes irr{\'e}ductibles de lieux sp{\'e}ciaux d'espaces de
					modules de courbes, action galoisienne en genre quelconque}},
			date={2015},
			ISSN={0373-0956},
			journal={Ann. Inst. Fourier (Grenoble)},
			volume={65},
			number={1},
			pages={245\ndash 276},
			url={http://aif.cedram.org/item?id=AIF_2015__65_1_245_0},
			review={\MR{3449153}},
		}
		
		\bib{COR87}{article}{
			author={Cornalba, Maurizio},
			title={{On the locus of curves with automorphisms}},
			date={1987},
			ISSN={0003-4622},
			journal={Ann. Mat. Pura Appl. (4)},
			volume={149},
			pages={135\ndash 151},
			url={http://dx.doi.org/10.1007/BF01773930},
			review={\MR{932781 (89b:14038)}},
		}
		
		\bib{Curran}{article}{
			author={Curran, M.~J.},
			title={{Automorphisms of semidirect products}},
			date={2008},
			ISSN={1393-7197},
			journal={Math. Proc. R. Ir. Acad.},
			volume={108},
			number={2},
			pages={205\ndash 210},
			url={http://dx.doi.org/10.3318/PRIA.2008.108.2.205},
			review={\MR{2475812}},
		}
		
		\bib{DD06}{article}{
			author={D{\`e}bes, Pierre},
			author={Emsalem, Michel},
			title={{Harbater-{M}umford components and towers of moduli spaces}},
			date={2006},
			ISSN={1474-7480},
			journal={J. Inst. Math. Jussieu},
			volume={5},
			number={3},
			pages={351\ndash 371},
			url={http://dx.doi.org/10.1017/S1474748006000053},
			review={\MR{2243760}},
		}
		
		\bib{DM69}{article}{
			author={Deligne, Pierre},
			author={Mumford, David},
			title={{The irreducibility of the space of curves of given genus}},
			date={1969},
			journal={Publications Math{\'e}matiques de l'IHES},
			volume={36},
			number={1},
			pages={75\ndash 109},
			review={\MR{262240}},
		}
		
		\bib{FRIED77}{article}{
			author={Fried, M.},
			title={{Fields of definition of function fields and {H}urwitz
					families---groups as {G}alois groups}},
			date={1977},
			ISSN={0092-7872},
			journal={Comm. Algebra},
			volume={5},
			number={1},
			pages={17\ndash 82},
			review={\MR{453746}},
		}
		
		\bib{FRI95}{incollection}{
			author={Fried, Michael~D.},
			title={{Introduction to modular towers: generalizing dihedral
					group--modular curve connections}},
			date={1995},
			booktitle={{Recent developments in the inverse {G}alois problem ({S}eattle,
					{WA}, 1993)}},
			series={{Contemp. Math.}},
			volume={186},
			publisher={Amer. Math. Soc., Providence, RI},
			pages={111\ndash 171},
			url={http://dx.doi.org/10.1090/conm/186/02179},
			review={\MR{1352270}},
		}
		
		\bib{TOH57}{article}{
			author={Grothendieck, Alexander},
			title={{Sur quelques points d'alg{\`e}bre homologique, {I \& II} }},
			date={1957},
			journal={Tohoku Mathematical Journal},
			volume={9},
			number={2-3},
			pages={119\ndash 221},
		}
		
		\bib{EGA4.2}{article}{
			author={Grothendieck, A.},
			title={{{{\'E}}l{\'e}ments de g{\'e}om{\'e}trie alg{\'e}brique. {IV}.
					{{\'E}}tude locale des sch{\'e}mas et des morphismes de sch{\'e}mas. {II}}},
			date={1965},
			ISSN={0073-8301},
			journal={Inst. Hautes {\'E}tudes Sci. Publ. Math.},
			number={24},
			pages={231},
		}
		
		\bib{EGA4.4}{article}{
			author={Grothendieck, Alexander},
			title={{{\'E}l{\'e}ments de g{\'e}om{\'e}trie alg{\'e}brique. {IV}.
					{\'E}tude locale des sch{\'e}mas et des morphismes de sch{\'e}mas {IV}}},
			date={1967},
			journal={Inst. Hautes {\'E}tudes Sci. Publ. Math.},
			number={32},
			review={\MR{0238860 (39 \#220)}},
		}
		
		\bib{ClassesChern}{incollection}{
			author={Grothendieck, A.},
			title={{Classes de {C}hern et repr{\'e}sentations lin{\'e}aires des
					groupes discrets}},
			date={1968},
			booktitle={{Dix expos{\'e}s sur la cohomologie des sch{\'e}mas}},
			series={{Adv. Stud. Pure Math.}},
			volume={3},
			publisher={North-Holland, Amsterdam},
			pages={215\ndash 305},
			review={\MR{265370}},
		}
		
		\bib{MaugeaisBordeaux}{article}{
			author={Maugeais, Sylvain},
			title={{Quelques calculs d'espaces {$\rm{R}^if_*G$} sur des courbes}},
			date={2016},
			ISSN={1246-7405},
			journal={Journal de Th{\'e}orie des Nombres de Bordeaux},
			volume={28},
			number={2},
			pages={361\ndash 390},
			url={http://jtnb.cedram.org/item?id=JTNB_2016__28_2_361_0},
			review={\MR{3509715}},
		}
		
		\bib{Milne}{book}{
			author={Milne, James~S.},
			title={{{\'E}tale cohomology}},
			series={{Princeton Mathematical Series}},
			publisher={Princeton University Press, Princeton, N.J.},
			date={1980},
			volume={33},
			ISBN={0-691-08238-3},
			review={\MR{559531}},
		}
		
		\bib{Rom09}{article}{
			author={Romagny, Matthieu},
			title={{Composantes connexes et irr{\'e}ductibles en familles}},
			date={2011},
			ISSN={0025-2611},
			journal={Manuscripta Math.},
			volume={136},
			number={1-2},
			pages={1\ndash 32},
			url={http://dx.doi.org/10.1007/s00229-010-0424-7},
			review={\MR{2820394}},
		}
		
		\bib{SER08}{book}{
			author={Serre, Jean-Pierre},
			title={{Topics in Galois Theory}},
			edition={{S}econd},
			series={{Research in Mathematics}},
			publisher={A K Peters LTD},
			date={2008},
			number={1},
		}
		
		\bib{Wells}{article}{
			author={Wells, Charles},
			title={Automorphisms of group extensions},
			date={1971},
			ISSN={0002-9947},
			journal={Trans. Amer. Math. Soc.},
			volume={155},
			pages={189\ndash 194},
			url={https://doi.org/10.2307/1995472},
			review={\MR{0272898}},
		}
		
		\bib{WEW98}{thesis}{
			author={Wewers, Stefan},
			title={{Construction of Hurwitz Spaces}},
			type={Ph.D. Thesis},
			address={Germany},
			date={1998},
		}
		
	\end{biblist}
\end{bibdiv}
\end{document}